\documentclass[a4paper,11pt,reqno]{amsart}
\usepackage{amsmath,amssymb}
\usepackage{mathrsfs}
\usepackage{amsfonts}
\usepackage[dvipsnames]{xcolor}
\usepackage{graphicx}
\usepackage{float}  
\usepackage{cite}
\usepackage{cases}
\usepackage{comment}
\usepackage{indentfirst}
\allowdisplaybreaks

\theoremstyle{plain}
\newtheorem{theorem}{Theorem}[section]
\newtheorem{lemma}[theorem]{Lemma}
\newtheorem{proposition}[theorem]{Proposition}
\newtheorem{corollary}[theorem]{Corollary}

\theoremstyle{definition}

\newtheorem{example}[theorem]{Example}

\theoremstyle{remark}

\numberwithin{equation}{section}

\usepackage[textwidth=16.0cm,textheight=21.0cm,centering]{geometry}
\usepackage[colorlinks=true,citecolor=blue,linkcolor=blue]{hyperref}

\begin{document}
\title[Interior regularity for the $k$-Hessian equation]
{Interior $W^{2,p}$ regularity for the $k$-Hessian equation}
\author[Y. Lyu]{Yasheng Lyu}
\address{School of Mathematics and Statistics, Xi'an Jiaotong University, Xi'an, Shaanxi 710049, People's Republic of China}
\email{lvysh21@stu.xjtu.edu.cn}

\date{}
\subjclass[2020]{Primary 35J60; Secondary 35B65, 35D30.}
\keywords{$k$-Hessian equation, Hessian measures, $W^{2,1}$ regularity, interior $W^{2,p}$ regularity, interior $C^{2,\alpha}$ regularity.}

\begin{abstract}
In this paper, we extend Caffarelli's interior $W^{2,p}$ and $C^{2,\alpha}$ regularity theory for the Monge-Amp\`ere equation to the $k$-Hessian equation, $3\leq k\leq n-1$, under sharp Sobolev or H\"older assumptions on the solution.
For densities bounded above and below by positive constants, we establish $W_{\mathrm{loc}}^{2,1}$ regularity for continuous $k$-convex solutions; the exponent $1$ is optimal even for constant density.
Higher interior regularity holds under either of two assumptions: $u\in W_{\mathrm{loc}}^{2,k(k-1)/2}$, or $u\in C_{\mathrm{loc}}^{1,\beta}$ for some $1-2/k<\beta<1$.
Under either assumption, if the density is positive and continuous, then $u\in W_{\mathrm{loc}}^{2,p}$ for every $1<p<\infty$; if the density is positive and belongs to $C_{\mathrm{loc}}^{\alpha}$, then $u\in C_{\mathrm{loc}}^{2,\alpha}$ for $0<\alpha<1$.
For densities bounded above and below by positive constants, we also establish global $W^{2,1}$ regularity on bounded strictly $(k-1)$-convex $C^{1,1}$ domains with $C^{1,1}$ boundary values.
\end{abstract}

\maketitle

\section{Introduction}

In this paper, we study the interior regularity of solutions to the $k$-Hessian equation
\begin{equation}\label{eqn1.1}
\sigma_{k}(D^{2}u)=f\quad\text{in }\Omega,
\end{equation}
where $\Omega\subset\mathbb{R}^{n}$ is a domain and $2\leq k\leq n$.
For $\lambda=(\lambda_{1},\dots,\lambda_{n})\in\mathbb{R}^{n}$, the $k$-th elementary symmetric polynomial is
\[
\sigma_{k}(\lambda):=\sum_{1\leq i_{1}<\cdots<i_{k}\leq n}\lambda_{i_{1}}\cdots\lambda_{i_{k}},\quad \forall 1\leq k\leq n.
\]
For a real symmetric matrix $A$ with eigenvalue vector $\lambda(A)$, we write $\sigma_{k}(A)=\sigma_{k}(\lambda(A))$.

The G\aa rding cones are given by
\[
\Gamma_{k}:=\{\lambda\in\mathbb{R}^{n}:\ \sigma_{i}(\lambda)>0,\ 1\leq i\leq k\},\quad \forall 1\leq k\leq n.
\]
Following Caffarelli-Nirenberg-Spruck \cite{CNS}, a $C^{2}$ function is $k$-admissible if $\lambda(D^{2}u)(x)\in\Gamma_{k}$ at every $x\in\Omega$. 
For $u\in C^{0}(\Omega)$ we use the weak notion of $k$-convexity: $u\in\Phi^{k}(\Omega)$ if, whenever a $C^{2}$ function $\varphi$ touches $u$ from above at $x\in\Omega$, one has $D^{2}\varphi(x)\in\overline{\Gamma}_{k}$. 
For such $u$, equation \eqref{eqn1.1} means
\[
\mu_{k}[u]=f\,dx\quad \text{in}\ \Omega,
\]
where $\mu_{k}[u]$ is the Hessian measure of Trudinger-Wang \cite{TW1,TW2}. 
If $u$ is $C^{2}$, then $\mu_{k}[u]=\sigma_{k}(D^{2}u)\,dx$. 
For further background, see Wang \cite{WangSurvey}.

For the Monge-Amp\`ere equation, Caffarelli \cite{CaffLocal,Caff} used compact sections of strictly convex solutions to prove interior $W^{2,p}$ regularity for all $1<p<\infty$ with positive continuous density, and interior $C^{2,\alpha}$ regularity with positive $C_{\mathrm{loc}}^{\alpha}$ density. 
We seek conditions that yield these results for general $k$-convex solutions when $3\leq k\leq n-1$.

For the Monge-Amp\`ere equation with the density bounded above and below by positive constants, De Philippis-Figalli \cite{DPF} established interior $W^{2,1}$ regularity of strictly convex Alexandrov solutions.
De Philippis-Figalli-Savin \cite{DPFS} and, independently, Schmidt \cite{Schmidt} improved the exponent to $1+\varepsilon$.
Without assuming strict convexity, Mooney proved interior $W^{2,1}$ regularity \cite{MooneyPartial} and subsequently obtained a local $L\log^{\delta}L$ estimate for $\Delta u$ with some $\delta=\delta(n)>0$ \cite{MooneyW21}.
For the $k$-Hessian equation with the density one, Hu \cite{Hu} proved interior $W^{2,1}$ regularity for convex viscosity solutions. 
For $k=2$, the author \cite{LyuSigma2Interior} obtained interior $W^{2,1}$ regularity for continuous $2$-convex solutions when the density is bounded above and below by positive constants.
Under the same density bounds, the author \cite{LyuSigma2Higher} also proved interior $W^{2,1+\varepsilon}$ regularity for some $\varepsilon>0$.
We extend the interior $W^{2,1}$ regularity result to continuous $k$-convex solutions for $3\leq k\leq n-1$.

\begin{theorem}\label{thm1.1}
Let $2\leq k\leq n$, and let $u\in C^{0}(\Omega)\cap\Phi^{k}(\Omega)$ satisfy
\begin{equation}\label{eqn1.2}
\mu_{k}[u]=f\,dx\quad\text{in }\Omega,\quad
0<\lambda\leq f\leq\Lambda<\infty\quad\text{a.e.\ in }\Omega.
\end{equation}
Then $u\in W_{\mathrm{loc}}^{2,1}(\Omega)$. 
Moreover, for every ball $B_{2r}\Subset\Omega$,
\[
\left\|D^{2}u\right\|_{L^{1}(B_{r})}\leq C(n)r^{n-2}\|u\|_{L^\infty(B_{2r})}.
\]
\end{theorem}

Mooney \cite[Sections 5 and 6]{MooneyW21} constructed a three-dimensional convex Alexandrov solution $v$ of $\det D^{2}v=1$ such that $v\notin W_{\mathrm{loc}}^{2,p}$ for any $p>1$.
For $3\leq k\leq n$, the extension $u(x,y,z)=v(x)+|y|^{2}/2$, with $y\in\mathbb{R}^{k-3}$ and $z\in\mathbb{R}^{n-k}$, satisfies $\mu_{k}[u]=dx$ and has the same failure of integrability, so the exponent $1$ in Theorem \ref{thm1.1} is sharp even for $f\equiv1$.

Mooney's example shows that bounds on the density alone do not yield higher integrability of the Hessian. 
To obtain higher interior regularity, we impose additional assumptions on $u$.
The critical Sobolev and H\"older exponents are
\begin{equation}\label{eqn1.3}
q_{*}:=\frac{k(k-1)}{2}\quad \text{and}\quad \beta_{*}:=1-\frac{2}{k}.
\end{equation} 
For positive $C^{1,1}$ densities, Urbas derived interior Hessian bounds assuming $u\in W_{\mathrm{loc}}^{2,q}(\Omega)$, first for $q>kn/2$ \cite{UrbasFirst} and later for $q>k(n-1)/2$ \cite{UrbasSecond}.
Under the latter condition, Urbas \cite{Urbas} also proved interior $C^{2,\alpha}$ regularity for every $0<\alpha<1$ when the density is positive and Lipschitz.
For positive $C^{1,1}$ densities, Urbas noted that an estimate of Liu-Bao \cite{LB} gives an interior Hessian bound when $q>k(k-1)$; for $k<n/2$, Dinew-Ko{\l}odziej \cite{DK} obtained such a bound when $q>n(k-1)/2$.
Urbas observed that, for $k\geq3$, no exponent $q<q_{*}$ can suffice for an interior Hessian bound and asked for the optimal exponent \cite[p.~221]{Urbas}.

In the real Monge-Amp\`ere case $k=n\geq3$, Urbas \cite{UrbasMA} proved that a positive lower bound on the density implies strict convexity if $u\in C_{\mathrm{loc}}^{1,\beta}$ with $\beta>1-2/n$ or if $u\in W_{\mathrm{loc}}^{2,q}$ with $q>n(n-1)/2$.
Collins-Mooney \cite{CM} obtained strict convexity at the Sobolev endpoint $q=n(n-1)/2$ under the same lower bound.

For $3\leq k\leq n-1$, we obtain higher interior regularity under the endpoint assumption $u\in W_{\mathrm{loc}}^{2,q_{*}}(\Omega)$ or, alternatively, under the assumption $u\in C_{\mathrm{loc}}^{1,\beta}(\Omega)$ for some $\beta_{*}<\beta<1$.
The key point is that, for positive continuous densities, either assumption gives strict separation: at the center of each ball $B\Subset\Omega$, $u$ lies strictly below the solution of the homogeneous $k$-Hessian equation in $B$ with boundary values $u$.
This separation allows us to construct a weight for the interior Hessian estimate, as in Chou-Wang \cite[Section~4]{CW}, and to localize the $W^{2,p}$ estimate in the following theorem.

\begin{theorem}\label{thm1.2}
Let $2\leq k\leq n$, and let $u\in C^{0}(\Omega)\cap\Phi^{k}(\Omega)$ satisfy
\begin{equation}\label{eqn1.4}
\mu_{k}[u]=f\,dx\quad\text{in }\Omega,\quad 0<\lambda\leq f\leq\Lambda<\infty, \quad f\in C^{0}(\Omega).
\end{equation} 
Assume that either
$u\in W_{\mathrm{loc}}^{2,q_{*}}(\Omega)$ or
$u\in C_{\mathrm{loc}}^{1,\beta}(\Omega)$ for some $\beta_{*}<\beta<1$.
Then $u\in W_{\mathrm{loc}}^{2,p}(\Omega)$ for every $1<p<\infty$.
For any such $p$ and any concentric balls $B_r\Subset B_R\Subset\Omega$,
\[
\|u\|_{W^{2,p}(B_r)}\leq C,
\]
where $C$ depends only on $n,k,\lambda,\Lambda,p,r,R$, $\|u\|_{L^\infty(B_R)}$, and a modulus of continuity of $f$ on $\overline{B_{R}}$, with the following additional dependence when $k\geq3$: on $\|u\|_{W^{2,q_{*}}(B_{R})}$ under the Sobolev assumption, or on $\beta$ and $[\nabla u]_{C^{0,\beta}(B_{R})}$ under the H\"older assumption.
\end{theorem}

When $k=2$, Theorem \ref{thm1.1} supplies the Sobolev hypothesis in Theorem \ref{thm1.2} since $q_{*}=1$. 
We thus recover the interior $W^{2,p}$ regularity of the author \cite{LyuSigma2Interior}. For the constant-density equation $\sigma_2(D^2u)=1$, Mooney \cite{MooneyQuadratic} showed that a $W^{2,p}$ bound with $p>2$ gives a pointwise interior Hessian bound. Warren-Yuan \cite{WY} obtained an interior Hessian estimate in dimension three, and Shankar-Yuan \cite{SYFour} established the four-dimensional estimate. 
Li-Wu \cite{LiWu} proved the estimate in every dimension.

For the quadratic Hessian equation with positive $C^{1,1}$ density, Qiu \cite{QiuThree} and Fan \cite{FanFour} proved interior Hessian estimates in dimensions three and four, respectively. 
In arbitrary dimension, Chen-Zhou-Zhu \cite{CZZ} established interior $C^{2,\alpha}$ regularity for admissible viscosity solutions with positive $C^\alpha$ density. 
For positive Dini continuous density, Hung-Liu-Trudinger \cite{HLT} obtained an interior Hessian bound for classical admissible solutions under a gradient bound.

For strictly convex solutions of the Monge-Amp\`ere equation, Caffarelli \cite{Caff} established interior $C^{2,\alpha}$ regularity for positive $C^\alpha$ density. 
Under the weaker assumption of positive Dini continuous density, Wang \cite{WangDini} proved that the Hessian is continuous, and Jian-Wang \cite{JW} estimated its modulus of continuity.

In the intermediate range $3\leq k\leq n-1$, we strengthen the continuity hypothesis on $f$ in Theorem \ref{thm1.2} to obtain an interior Hessian bound. 
The condition we use is the Dini condition on the square root of the modulus of $f$.
For a nondecreasing modulus of continuity $\omega$, set
\begin{equation}\label{eqn1.5}
\Psi_{\omega}(\rho):=\int_{0}^{\rho}\frac{\sqrt{\omega(r)}}{r}\,dr,\quad 0<\rho\leq 1.
\end{equation}
For each nonempty compact set $K\subset\Omega$, let
\[
\omega_{f,K}(r):=\sup\{|f(x)-f(y)|:\ x,y\in K,\ |x-y|\leq r\}.
\]
We say that $f$ is locally square-root Dini continuous if
$\Psi_{\omega_{f,K}}(1)<\infty$ for every such $K$.

\begin{theorem}\label{thm1.3}
Under the hypotheses of Theorem \ref{thm1.2}, assume in addition that
$f$ is locally square-root Dini continuous. 
Then
$u\in W_{\mathrm{loc}}^{2,\infty}(\Omega)$. 
For any concentric balls $B_{r}\Subset B_{R}\Subset\Omega$,
\[
\|u\|_{W^{2,\infty}(B_{r})}\leq C,
\]
where $C$ depends only on $n,k,\lambda,\Lambda,r,R$,
$\|u\|_{L^\infty(B_{R})}$, and the square-root Dini modulus of $f$ on $\overline{B_{R}}$, with the following additional dependence when $k\geq3$: on $\|u\|_{W^{2,q_{*}}(B_{R})}$ under the Sobolev assumption, or on $\beta$ and $[\nabla u]_{C^{0,\beta}(B_{R})}$ under the H\"older assumption.
\end{theorem}

\begin{corollary}\label{thm1.4}
Under the hypotheses of Theorem \ref{thm1.2}, suppose that $f\in C_{\mathrm{loc}}^{0,\alpha}(\Omega)$ for some $0<\alpha<1$.
Then $u\in C_{\mathrm{loc}}^{2,\alpha}(\Omega)$.
\end{corollary}

When $k=2$, we recover the interior $C^{2,\alpha}$ regularity proved by Chen-Zhou-Zhu \cite{CZZ}.

We also obtain global $W^{2,1}$ regularity for the Dirichlet problem. 
A $C^{1,1}$ boundary is strictly $(k-1)$-convex if its principal curvatures $\kappa$ satisfy $\sigma_{j}(\kappa)\geq c_{0}>0$ almost everywhere for $1\leq j\leq k-1$.
This is strict mean convexity for $k=2$ and uniform convexity for $k=n$.

\begin{theorem}\label{thm1.5}
Let $\Omega$ be bounded with strictly $(k-1)$-convex $C^{1,1}$ boundary. If $u\in C^{0}(\overline{\Omega})\cap\Phi^{k}(\Omega)$ satisfies \eqref{eqn1.2} with trace $\varphi\in C^{1,1}(\partial\Omega)$, then $u\in W^{2,1}(\Omega)$ and
\[
\|u\|_{W^{2,1}(\Omega)}\leq C\big(\|\varphi\|_{C^{1,1}(\partial\Omega)}+\Lambda^{1/k}\big).
\]
Here $C$ depends only on $n,k$, the diameter of $\Omega$, the radii and $C^{1,1}$ bounds of a boundary charts, and $c_{0}$.
\end{theorem}

The rest of the paper is organized as follows.
Section~2 proves Theorem~\ref{thm1.1}.
Section~3 establishes strict separation from the homogeneous solution.
Section~4 proves Theorem~\ref{thm1.2}.
Section~5 proves Theorem~\ref{thm1.3} and Corollary~\ref{thm1.4}.
Section~6 proves Theorem~\ref{thm1.5} and gives the boundary counterexample.
Appendix~A proves the solvability of the auxiliary Dirichlet problem.

\section{Interior $W^{2,1}$ regularity}

To prove Theorem \ref{thm1.1}, we show that $\Delta u$ is absolutely continuous. We compare the singular parts of $\Delta u$ and a mixed Hessian measure.

\subsection{The distributional Hessian}

Throughout, $G(A):=\sigma_{k}(A)^{1/k}$ on $\overline{\Gamma}_{k}$. It is concave, homogeneous of degree one, and elliptic in $\Gamma_{k}$ (see Caffarelli-Nirenberg-Spruck \cite{CNS}). We use the Frobenius norm $\|\cdot\|$ on symmetric matrices and the corresponding total variation $|M|$ for matrix-valued measures. The superscript $s$ denotes the singular part with respect to Lebesgue measure.

\begin{lemma}\label{thm2.1}
Let $2\leq k\leq n$ and $w\in C^{0}(\Omega)\cap\Phi^{k}(\Omega)$.
Then $w\in W_{\mathrm{loc}}^{1,1}(\Omega)$, and $D^{2}w$ is a matrix-valued Radon measure with
\begin{equation}\label{eqn2.1}
|D^{2}w|\leq\Delta w.
\end{equation}
Its density $S=dD^{2}w/d(\Delta w)$ satisfies $S\in\overline{\Gamma}_{k}$ and $\operatorname{tr} S=1$ $\Delta w$-a.e.
\end{lemma}
\begin{proof}
The gradient of $w$ is locally integrable by Trudinger-Wang \cite{TW1,TW2}. Subharmonicity gives the nonnegative Radon measure $\mu=\Delta w$. The Hessian of each radial mollification $w_{\varepsilon}$ lies in $\overline{\Gamma}_{k}\subset\overline{\Gamma}_{2}$, so
\[
\|D^{2}w_{\varepsilon}\|\leq\Delta w_{\varepsilon}
\]
by $2\sigma_{2}(H)=(\operatorname{tr}H)^{2}-\|H\|^{2}$. Passing to the limit gives \eqref{eqn2.1} and $D^{2}w=S\mu$. Taking the trace gives $\operatorname{tr}S=1$ $\mu$-a.e. For each linear functional nonnegative on $\overline{\Gamma}_{k}$, its evaluation on $D^{2}w$ is a nonnegative measure by passage to the limit. A countable dense family of such functionals gives $S\in\overline{\Gamma}_{k}$ $\mu$-a.e.
\end{proof}

\begin{lemma}\label{thm2.2}
Let $2\leq k\leq n$ and $u\in C^{0}(\Omega)\cap\Phi^{k}(\Omega)$. If $\mu_{k}[u]=f\,dx$ with $f\in L^{1}_{\mathrm{loc}}(\Omega)$ and $D^{2}u=A\,dx+(D^{2}u)^{s}$, then
\begin{equation}\label{eqn2.2}
\sigma_{k}(A)=f\quad\text{a.e.\ in }\Omega.
\end{equation}
\end{lemma}
\begin{proof}
For $k>n/2$, see Chaudhuri-Trudinger \cite[Theorem 3.2]{CT}. In the present setting we use an $L^{1}$ blow-up. For a.e.\ $x$, $A$ and $f$ have Lebesgue values at $x$, $r^{-n}|(D^{2}u)^{s}|(B_{r}(x))\to0$, and $A(x)\in\overline{\Gamma}_{k}$ by Lemma \ref{thm2.1}.
On $B_{1}$, define
\[
q(y):=\frac{1}{2} y^{T} A(x)y,\quad
w_{r}(y):=r^{-2}u(x+ry)-q(y).
\]
Then
\[
|D^{2}w_{r}|(B_{1})
=r^{-n}|D^{2}u-A(x)\, dx|(B_{r}(x))\rightarrow0.
\]
Subtract an affine function $\ell_{r}$ so that $w_{r}-\ell_{r}$ and its gradient have mean zero. The Poincar\'e inequalities in $W^{1,1}$ and $BV$ (see Ambrosio-Fusco-Pallara \cite[Theorem 3.44 and Remark 3.45]{AFP}) give
\[
\|w_{r}-\ell_{r}\|_{L^{1}(B_{1})}
\leq C\|\nabla w_{r}-\nabla\ell_{r}\|_{L^{1}(B_{1})}
\leq C|D^{2}w_{r}|(B_{1}).
\]
Thus $r^{-2}u(x+r\,\cdot)-\ell_{r}\to q$ in $L^{1}(B_{1})$. By scaling and affine invariance its Hessian measure is $f(x+r\,\cdot)\,dy$, as follows by local mollification. This tends to $f(x)\,dy$ by differentiation and to $\sigma_{k}(A(x))\,dy$ by weak continuity (see Trudinger-Wang \cite[Theorem 1.1(ii)]{TW2}).
\end{proof}

\subsection{Locally constant maximizing radii}

Let $u$ satisfy the hypotheses of Theorem \ref{thm1.1}, and write $D^{2}u=A\,dx+(D^{2}u)^{s}$. Fix $B_{4R}(x_{0})\Subset\Omega$, abbreviate $B_{s}=B_{s}(x_{0})$, and set
\[
D:=B_{2R},\quad \mu:=\Delta u=T\, dx+\mu^{s},\quad T:=\operatorname{tr} A.
\]
On $\overline{B}_{3R}$, $u$ is Lipschitz with constant $L$ by the gradient estimate for bounded $f$ (see Liu-Ma-Zang \cite[Corollary 1.4]{LMZ}; see also He \cite[Theorem 1.1]{He}), or by convexity when $k=n$.
Let $\rho\in C_{c}^{\infty}(B_{1}(0))$ be a nonnegative radial mollifier positive on $B_{7/8}(0)$, and set $u_{r}=u*\rho_{r}$ with $\rho_{r}(y)=r^{-n}\rho(y/r)$. Let $\mathcal{M}\mu$ be the centered Hardy--Littlewood maximal function of the measure $\mu|_{B_{3R}}$ on $\mathbb R^{n}$.
For $x\in D$ and $0<2r<R$ we have
\begin{equation}\label{eqn2.3}
cr^{2-n}\mu(B_{r/4}(x))
\leq u_{r}(x)-u(x)\leq Cr^{2}\mathcal{M}\mu(x),
\end{equation}
\begin{equation}\label{eqn2.4}
u_{2r}(x)-u_{r}(x)
\geq cr^{2-n}\mu(B_{r}(x)).
\end{equation}
Indeed, the spherical mean $m_{s}u(x)$ satisfies
\[
m_{s}u(x)-u(x)
=\frac{1}{|S^{n-1}|}\int_{0}^{s} t^{1-n}\mu(B_{t}(x))\, dt.
\]
The divergence theorem gives the formula between two positive radii by mollification, since ball masses converge at a.e.\ radius. Let the smaller radius tend to zero. Averaging $m_{ar}u$ against $|S^{n-1}|a^{n-1}\rho(ae_{1})$ gives $u_{r}$ and the maximal-function bound. For the lower bounds, restrict to $3/4<a<7/8$ and integrate over $[r/4,r/2]$ and $[r,3r/2]$, respectively.

Fix $0<r_{0}<R/4$ and let $r_{j}:=2^{-j}r_{0}$ for $j\geq0$.
For $P>L/r_{0}$ put
\[
v_{P}:=\max\left\{u,\ \sup_{j\geq0}(u_{r_{j}}-Pr_{j}^{2})\right\},
\quad E_{P}:=\{x\in D:v_{P}(x)>u(x)\}.
\]
The finite maxima with $u$ converge uniformly to $v_{P}$, so $v_{P}$ is continuous and $k$-convex (see Trudinger-Wang \cite[Lemmas 2.3 and 2.5]{TW2}). We have $u_{r}-u-Pr^{2}\leq Lr-Pr^{2}$ and, by \eqref{eqn2.3}, $E_{P}\subset\{\mathcal{M}\mu\geq cP\}$. Maximizing the quadratic and using the weak $(1,1)$ bound gives
\begin{equation}\label{eqn2.5}
0\leq v_{P}-u\leq\frac{L^{2}}{4P},\quad |E_{P}|\leq\frac{C}{P}.
\end{equation}
The constants below are independent of $P$.

At $\mu^{s}$-a.e.\ $x$, $s^{-n}\mu(B_{s}(x))\to\infty$ as $s\downarrow0$, by differentiation of Lebesgue measure with respect to $\mu+dx$ (see Ambrosio-Fusco-Pallara \cite[Theorem 2.22]{AFP}). The lower bound in \eqref{eqn2.3} then gives $\mu^{s}(D\setminus E_{P})=0$. For each $i\ne j$, the ratio $(u_{r_{i}}-u_{r_{j}})/(r_{i}^{2}-r_{j}^{2})$ has at most countably many level sets of positive $\mu^{s}$-measure. Take $P\to\infty$ outside these levels for all pairs.

For each such $P$ and $\mu^{s}$-a.e.\ $x$, the maximum comes from a unique radius $r=r_{j}$. Since $v_{P}(x)>u(x)$ and the terms converge uniformly to $u$, this remains the maximizing radius in a neighborhood of $x$. Hence $r$ is measurable there. The Lipschitz bound gives $r<L/P<r_{0}$, so the radius $2r$ occurs in the supremum. Comparing with its term gives $u_{2r}(x)-u_{r}(x)\leq3Pr^{2}$. By \eqref{eqn2.4},
\begin{equation}\label{eqn2.6}
\mu(B_{r}(x))\leq CPr^{n}
\end{equation}
for $\mu^{s}$-almost every $x\in D$ and every $P$ in the chosen sequence.

\subsection{Absolute continuity}

Write $\widetilde{\sigma}_{k}$ for the symmetric $k$-linear polarization of $\sigma_{k}$, normalized by $\widetilde{\sigma}_{k}(A,\ldots,A)=\sigma_{k}(A)$, and $\widetilde{\mu}_{k}$ for the corresponding mixed Hessian measure of Trudinger-Wang \cite{TW3}. With one entry equal to $u$, set
\[
\nu_{P}:=\widetilde{\mu}_{k}[u,v_{P},\ldots,v_{P}].
\]
We shall prove
\begin{equation}\label{eqn2.7}
(\nu_{P})^{s}(B_{R})\leq C/P.
\end{equation}
Choose $0\leq\chi\leq1$ in $C_{c}^{\infty}(D)$ equal to one on $B_{R}$.
We first prove
\begin{equation}\label{eqn2.8}
\int_{D}\chi\,d(\nu_{P}-\mu_{k}[u])\leq C/P.
\end{equation}
Mollify $u$ and $v_{P}$ near $K=\operatorname{supp}\chi$. Replace the last $k-1$ entries of the smooth mixed Hessian one at a time by $v_{P,\varepsilon}$. In each difference, the coefficient of $D^{2}(v_{P,\varepsilon}-u_{\varepsilon})$ is nonnegative definite and divergence free. Its trace is $(n-k+1)/k$ times the mixed $(k-1)$-Hessian of the remaining entries (see Trudinger-Wang \cite[(2.10)--(2.14)]{TW3}). Integrating twice by parts and setting $w_{\varepsilon}=u_{\varepsilon}+v_{P,\varepsilon}$ bounds the absolute value of the smooth integral by
\[
C\|D^{2}\chi\|_{L^{\infty}}
\int_{K}(v_{P,\varepsilon}-u_{\varepsilon})
\sigma_{k-1}(D^{2}w_{\varepsilon})\, dx.
\]
Here we used positivity of all terms in the mixed expansion of $\sigma_{k-1}(D^{2}w_{\varepsilon})$.

By \eqref{eqn2.5}, $0\leq v_{P,\varepsilon}-u_{\varepsilon}\leq L^{2}/(4P)$ on $K$ for small $\varepsilon$. The functions $w_{\varepsilon}$ are $(k-1)$-convex and uniformly bounded on a fixed neighborhood of $K$. The local mass bound (see Trudinger-Wang \cite[(2.12)]{TW1}) therefore gives
\[
\int_{K}\sigma_{k-1}(D^{2}w_{\varepsilon})\, dx\leq C
\]
uniformly in $\varepsilon$ and $P$. Letting $\varepsilon\downarrow0$ for fixed $P$ gives \eqref{eqn2.8} by weak continuity of mixed Hessian measures (see Trudinger-Wang \cite[Theorem 2.4]{TW3}).

On $D\setminus E_{P}$, the absolutely continuous Hessian densities of $v_{P}$ and $u$ agree a.e. by locality of approximate gradients, applied first to $v_{P}-u$ and then to the components of its gradient (see Ambrosio-Fusco-Pallara \cite[Proposition 3.73 and Theorem 3.83]{AFP}).
If $B$ denotes the absolutely continuous Hessian density of $v_{P}$, simultaneous mollification and Fatou's lemma, as in the quadratic case (see the author \cite[Lemma 2.1, (2.5)]{LyuSigma2Interior}), give
\[
\nu_{P}\geq\widetilde{\sigma}_{k}(A,B,\ldots,B)\, dx.
\]
Here the smooth polarized densities are nonnegative, converge almost everywhere to the displayed density, and converge weakly as measures to $\nu_{P}$.
Thus on $D\setminus E_{P}$ the density of $\nu_{P}$ is at least $\sigma_{k}(A)=f$ almost everywhere, by Lemma \ref{thm2.2}.
It follows from \eqref{eqn2.5} and \eqref{eqn2.8} that
\[
(\nu_{P})^{s}(B_{R})
\leq\int_{D}\chi\,d(\nu_{P}-\mu_{k}[u])+\int_{E_{P}}\chi f\, dx
\leq C/P.
\]

On an open set where $v_{P}=u_{r}-Pr^{2}$, the polarized form (see Trudinger-Wang \cite[(2.8), (2.11) and Theorem 2.4]{TW3}) gives
\begin{equation}\label{eqn2.9}
(\nu_{P})^{s}=\frac{1}{k}\sum_{i,j=1}^{n}
\sigma_{k}^{ij}(D^{2}u_{r})(D_{ij}u)^{s}.
\end{equation}
Here $D\sigma_{k}(H)[Z]=\sum_{i,j}\sigma_{k}^{ij}(H)Z_{ij}$ for symmetric $Z$.
The formula follows by mollifying $u$ in the first entry and testing against compactly supported functions in that open set.

\begin{proof}[Proof of Theorem \ref{thm1.1}]
By Lemma \ref{thm2.1}, write $(D^{2}u)^{s}=S\mu^{s}$, where $S\in\overline{\Gamma}_{k}$ and $\operatorname{tr} S=1$ $\mu^{s}$-almost everywhere.
Fix such a point $x$ with a unique maximizing radius $r$, and put
\[
B:=B_{r/2}(x),\quad
A_{B}:=\frac{1}{|B|}\int_{B} A(y)\, dy.
\]
By Lemma \ref{thm2.2} and the concavity of $G$,
\[
G(A_{B})\geq\frac{1}{|B|}\int_{B} G(A(y))\, dy\geq\lambda^{1/k}.
\]
Thus $A_{B}\in\Gamma_{k}$ and $\operatorname{tr} A_{B}=|B|^{-1}\int_{B} T\, dy$.
The lower bound for $\rho_{r}$ on $B$, together with the decomposition $D^{2}u=S\mu$ in Lemma \ref{thm2.1}, gives
\[
D^{2}u_{r}(x)=cA_{B}+H_{0},\quad H_{0}\in\overline{\Gamma}_{k},
\]
for a fixed $c=c(n,\rho)>0$.
Indeed, subtract $c/|B|$ from the convolution weight of the absolutely continuous part on $B$. The remaining weights are nonnegative and the matrix densities belong to $\overline{\Gamma}_{k}$.

For $H,Z\in\overline{\Gamma}_{k}$ we use
\begin{equation}\label{eqn2.10}
(\operatorname{tr}H)D\sigma_{k}(H)[Z]\geq\sigma_{k}(H)\operatorname{tr}Z.
\end{equation}
Indeed, $Q:=(\sigma_{k}/\sigma_{1})^{1/(k-1)}$ is concave and homogeneous of degree one. For $H,Z\in\Gamma_{k}$ one has $DQ(H)[Z]\geq Q(Z)>0$. Logarithmic differentiation of $Q$ gives \eqref{eqn2.10}, first in the open cone and then on its closure by continuity.

Expanding in $H_{0}$ and using positivity of $\widetilde{\sigma}_{k}$ on $\overline{\Gamma}_{k}$, we obtain
\begin{equation}\label{eqn2.11}
D\sigma_{k}(D^{2}u_{r}(x))[S(x)]
\geq c^{k-1}D\sigma_{k}(A_{B})[S(x)]
\geq c^{k-1}\frac{\sigma_{k}(A_{B})}{\operatorname{tr} A_{B}}
\geq c^{k-1}\lambda\frac{|B|}{\int_{B} T\, dy}.
\end{equation}
The second inequality is \eqref{eqn2.10} with $\operatorname{tr} S(x)=1$.

For each $P$ the sets on which the maximizing radius is locally constant cover $\mu^{s}$-almost all of $B_R$. By \eqref{eqn2.7}, \eqref{eqn2.9} and \eqref{eqn2.11},
\[
\int_{B_{R}}
\frac{P|B_{r/2}(x)|}{\int_{B_{r/2}(x)}T\, dy}\, d\mu^{s}(x)
\leq CP(\nu_{P})^{s}(B_{R})\leq C.
\]
On a common set of full $\mu^{s}$-measure, $r=r(x,P)<L/P$. Differentiation of $T\,dx$ with respect to $\mu$ gives
\[
\frac{\int_{B_{r/2}(x)}T\, dy}{\mu(B_{r/2}(x))}
\rightarrow0
\quad\text{as }P\rightarrow\infty,\quad \mu^{s}\text{-a.e.\ }x.
\]
Since $\mu(B_{r/2}(x))\leq CPr^{n}$ by \eqref{eqn2.6}, the integrand tends to infinity $\mu^{s}$-a.e. Fatou's lemma gives $\mu^{s}(B_{R})=0$. Varying the ball and using \eqref{eqn2.1} yields $u\in W_{\mathrm{loc}}^{2,1}(\Omega)$.

For the norm estimate, write $B_{s}=B_{s}(x_{0})$ and put $M=\|u\|_{L^{\infty}(B_{2r})}$.
Choose $\eta\in C^{\infty}_{c}(B_{2r})$ with $0\leq\eta\leq1$, $\eta=1$ on $B_{r}$ and $\|D^{2}\eta\|_{L^{\infty}}\leq C(n)r^{-2}$. By \eqref{eqn2.1},
\[
\int_{B_{r}}\|D^{2}u\|\,dx
\leq\int\eta\,d\Delta u
=\int u\Delta\eta\,dx
\leq C(n)r^{n-2}M.
\]
Let $\zeta\in C^{\infty}_{c}(B_{r})$ be a rescaling of a fixed nonnegative bump of integral one on the unit ball. Then $|\int\zeta\nabla u|=|\int u\nabla\zeta|\leq C(n)r^{-1}M$.
The Poincar\'e inequality for $\nabla u$ with this weighted mean gives
\[
\|\nabla u\|_{L^{1}(B_{r})}
\leq C(n)r\|D^{2}u\|_{L^{1}(B_{r})}+|B_{r}|\Big|\int\zeta\nabla u\Big|
\leq C(n)r^{n-1}M.
\]
Together with $\|u\|_{L^{1}(B_{r})}\leq |B_{r}|M$, this proves the estimate.
\end{proof}

\section{Strict separation from the homogeneous solution}

Let $2\leq k\leq n$ and $u\in C^{0}(\Omega)\cap\Phi^{k}(\Omega)$ satisfy \eqref{eqn1.4}. For a ball $B=B_{R}(x_{0})\Subset\Omega$, let $H_{B}u$ denote the solution of the homogeneous Dirichlet problem with boundary values $u$. We prove $H_{B}u(x_{0})>u(x_{0})$ under either regularity condition in Theorem \ref{thm1.2}. This gap supplies the weight in Section 4.

\subsection{The homogeneous Dirichlet problem}

By Trudinger-Wang \cite[Theorem 3.5]{TW2}, there is a unique continuous $k$-convex solution $H_{B}u$ of
\begin{equation}\label{eqn3.2}
\mu_{k}[H_{B}u]=0\quad\text{in }B,\qquad H_{B}u=u\quad\text{on }\partial B.
\end{equation}
By comparison, it is the largest continuous $k$-convex function on $\overline B$ with boundary values at most $u$. Thus $u\leq H_{B}u$, and $H_{B'}u\leq H_{B}u$ on $B'\Subset B$. For $k=n$, $H_{B}u$ is the convex envelope of the boundary values.

If $Q\in C^{2}$ touches $H_{B}u$ from below at $y\in B$, then $D^{2}Q(y)\notin\Gamma_{k}$. Indeed, if $Q$ were admissible near $y$, adding a small positive multiple of $s^{2}/2-|x-y|^{2}$ would contradict comparison on a small ball $B_{s}(y)$.

For smooth admissible $v_{i}$ with $G(D^{2}v_{i})=h_{i}$, comparison gives
\begin{equation}\label{eqn3.1}
\|v_{1}-v_{2}\|_{L^{\infty}(B)}\leq \|v_{1}-v_{2}\|_{L^{\infty}(\partial B)}
+\frac{R^{2}}{2G(I)}\|h_{1}-h_{2}\|_{L^{\infty}(B)}.
\end{equation}
Indeed, put $d=\|h_{1}-h_{2}\|_{L^{\infty}(B)}$ and $q=d(|x-x_{0}|^{2}-R^{2})/(2G(I))$. Concavity and homogeneity give $G(D^{2}(v_{1}+q))\geq h_{1}+d\geq h_{2}$. Compare $v_{1}+q$ with $v_{2}+\|v_{1}-v_{2}\|_{L^{\infty}(\partial B)}$ and interchange the indices.

For a nonempty compact set $K\subset B$, $u<H_{B}u$ on $K$ is equivalent to the existence of a smooth admissible $W$ on $\overline B$ with
\begin{equation}\label{eqn3.3}
W>u\quad\text{on }K,\qquad W<u\quad\text{on }\partial B.
\end{equation}
Comparison gives $W\leq H_{B}u$. Conversely, let $\varphi_{j}$ be smooth boundary data converging uniformly to $u$ on $\partial B$. By Caffarelli-Nirenberg-Spruck \cite{CNS}, there are smooth admissible solutions $v_{j}$ with $G(D^{2}v_{j})=1/j$ in $B$ and $v_{j}=\varphi_{j}$ on $\partial B$. Estimate \eqref{eqn3.1} gives uniform convergence, and weak continuity and uniqueness identify the limit with $H_{B}u$ (see Trudinger-Wang \cite{TW1,TW2}). Take $W=v_{j}-2s$ with $0<4s<\min_{K}(H_{B}u-u)$ and $\|v_{j}-H_{B}u\|_{L^{\infty}(B)}<s$. The two strict inequalities persist under small uniform perturbations of $u$. For convex solutions of $\sigma_{2}(D^{2}u)=1$, Mooney \cite[Section 2]{MooneySigma2} constructed a quadratic comparison function satisfying \eqref{eqn3.3} on an inner ball.

Separation on all sufficiently small centered balls excludes harmonic restrictions to open disks in affine planes of dimension $n-k+1$. Indeed, for $H\in\Gamma_{j}$ and a unit vector $e$, ellipticity gives
\[
\sigma_{i}(H|_{e^{\perp}})=D\sigma_{i+1}(H)[ee^{T}]>0,
\quad 1\leq i\leq j-1.
\]
Restricting $D^{2}W$ successively thus makes $W$ subharmonic on the intersection of the plane with a ball centered in the disk. If $u$ were harmonic there, \eqref{eqn3.3} at the center would contradict the maximum principle on that intersection. For convex $u$, the contact set of any supporting affine function therefore has affine dimension at most $n-k$. Mooney proved this bound for $k=2$ \cite[Theorem 1.1]{MooneySigma2} and, for $\sigma_{2}=1$, excluded contact from below with a harmonic function on a smooth open portion of a minimal hypersurface \cite[Theorem 1.1]{MooneyQuadratic}.

For $k=n$, the dimension bound excludes affine segments. Conversely, if $u$ is strictly convex on a ball, a supporting plane at its center lies strictly below $u$ on the sphere. Adding a small positive constant and then a smaller positive quadratic gives \eqref{eqn3.3} at the center. For $k<n$, however, strict separation does not require convexity of $u$.

\subsection{The H\"older criterion}

Put $g=f^{1/k}$. Let $g_{j}>0$ and $\varphi_{j}$ be smooth uniform approximations of $g$ on $\overline B$ and $u$ on $\partial B$. By Caffarelli-Nirenberg-Spruck \cite{CNS}, there is a smooth admissible $U_{j}$ with $G(D^{2}U_{j})=g_{j}$ and boundary value $\varphi_{j}$. By \eqref{eqn3.1}, $U_{j}$ converges uniformly, and weak continuity and comparison identify the limit with $u$ (see Trudinger-Wang \cite{TW1,TW2}).

Let $U_{\varepsilon}=u*\zeta_{\varepsilon}$, where $\zeta_{\varepsilon}(x)=\varepsilon^{-n}\zeta(x/\varepsilon)$ for a nonnegative radial mollifier $\zeta\in C_{c}^{\infty}(B_{1})$. Concavity gives $G(D^{2}(U_{j}*\zeta_{\varepsilon}))\geq g_{j}*\zeta_{\varepsilon}$. The convolutions converge locally in $C^{2}$ for fixed $\varepsilon$, so $G(D^{2}U_{\varepsilon})\geq g*\zeta_{\varepsilon}$. If $f\geq\lambda>0$ on $\overline B$, this and subharmonicity give
\begin{equation}\label{eqn3.4}
D^{2}U_{\varepsilon}\in\Gamma_{k},\quad
\sigma_{k}(D^{2}U_{\varepsilon})\geq\lambda,\quad U_{\varepsilon}\geq u,
\end{equation}
wherever $\overline B_{\varepsilon}(x)\subset B$.

\begin{lemma}\label{thm3.1}
Let $2\leq k\leq n$ and let $u\in \Phi^{k}(\Omega)$ satisfy
$\mu_{k}[u]=f\,dx$ in $\Omega$, where $f\in C^{0}(\Omega)$ and $f>0$.
Suppose that $u\in C_{\mathrm{loc}}^{1,\beta}(\Omega)$ for some $1-2/k<\beta<1$.
Then $H_{B}u(x_{0})>u(x_{0})$ for every ball $B=B_{R}(x_{0})\Subset\Omega$.
\end{lemma}

\begin{proof}
Mollify on a slightly larger ball and put $\delta=\varepsilon^{2-2/k}$. Cancellation of the linear part of $u$ under radial averaging and the $C^{0,\beta}$ bound for $Du$ give
\[
\delta^{-1}\|U_{\varepsilon}-u\|_{L^{\infty}(B)}
=O(\varepsilon^{\beta-\beta_{*}}),\quad
\delta\|D^{2}U_{\varepsilon}\|_{L^{\infty}(B)}^{k-1}
=O(\varepsilon^{(k-1)(\beta-\beta_{*})}).
\]
Both quantities tend to zero. Since $D\sigma_{k}$ is a polynomial of degree $k-1$, subtracting $2t\delta I$ from $D^{2}U_{\varepsilon}$ changes $\sigma_{k}$ by $o(1)$ uniformly for $0\leq t\leq1$. By \eqref{eqn3.4}, $\sigma_{k}$ remains positive along this segment. It stays in $\Gamma_{k}$ since Maclaurin's inequalities give $\sigma_{k}=0$ on $\partial\Gamma_{k}$. The function $W=U_{\varepsilon}+\delta(R^{2}/2-|x-x_{0}|^{2})$ is admissible and satisfies \eqref{eqn3.3} with $K=\{x_{0}\}$.
\end{proof}

\subsection{Barriers and mixed norms}

At the Sobolev endpoint, we use mixed-norm bounds for the mollifications and a convex barrier.

Put $d=k-1$, $m=n-k+1$ and $q=kd/2$. For each set $I$ of $m$ coordinate indices, write $x=(z,y)$ with $z=\pi_{I}x=x_{I}$ and $y=x_{I^{c}}$. Write $\mathcal{K}_{k}$ for the trace-one section of the dual cone of $\overline{\Gamma}_{k}$ and set $L_{A}w:=\operatorname{tr}(AD^{2}w)$. Every $A\in\mathcal{K}_{k}$ is positive semidefinite since $\overline{\Gamma}_{k}$ contains the positive semidefinite cone. The set $\mathcal{K}_{k}$ is compact, and
\[
H\in\Gamma_{k}\quad\Longleftrightarrow\quad
\operatorname{tr}(AH)>0\quad\text{for every }A\in\mathcal{K}_{k}.
\]
Indeed, the corresponding nonstrict inequalities characterize the closed cone by duality (see Trudinger-Wang \cite[Lemma 2.2]{TW2}). Applying this characterization to $H-\eta I$ for $\eta>0$ and using compactness gives the equivalence. For maximum estimates associated with the dual G\aa rding cone, see Kuo-Trudinger \cite[Theorem 1.1 and Proposition 2.1]{KT}.

\begin{lemma}\label{thm3.2}
For $0\leq g\in C_{c}^{0}(B_{1/2})$, set $g_{I}(z)=\sup_{y}g(z,y)$ using the zero extension. There is a smooth convex $J$ on $B_{1}$ with
\[
-C(n,k)\sum_{I}\|g_{I}\|_{L^{m}}\leq J\leq0\quad\text{in }B_{1},
\qquad L_{A}J\geq g\quad\text{in }B_{1/2}\quad(A\in\mathcal{K}_{k}).
\]
\end{lemma}

\begin{proof}
Each $A\in\mathcal{K}_{k}$ has rank at least $m$. Otherwise its kernel contains a $k$-plane with orthogonal projection $P$. Since $P\in\Gamma_{k}$, we have $P-sI\in\Gamma_{k}$ for small $s>0$, whereas $\operatorname{tr}(A(P-sI))=-s$. Thus $\sigma_{m}(A)>0$ on the compact set $\mathcal{K}_{k}$. As $\sigma_{m}(A)$ is the sum of the principal $m$-minors, each $A$ has a principal block $A_{II}$ with $\det A_{II}\geq c(n,k)>0$.

If $g=0$, take $J=0$. Otherwise each $g_{I}$ is continuous and nonzero. On the unit $m$-ball choose a smooth positive majorant $f_{I}\geq g_{I}$ with $\|f_{I}\|_{L^{m}}\leq2\|g_{I}\|_{L^{m}}$. Let $J_{I}$ be the convex solution of $\det D^{2}J_{I}=f_{I}^{m}$ with zero boundary values. For $m=1$ this is $J_{I}''=f_{I}$. Then $J_{I}\leq0$ and $\|J_{I}\|_{L^{\infty}}\leq C_{m}\|f_{I}\|_{L^{m}}$, since the gradient image contains the ball of slopes of length less than $\|J_{I}\|_{L^{\infty}}/2$ and has measure at most $\int f_{I}^{m}$. For each $A$, the block selected above satisfies
\[
L_{A}(J_{I}\circ\pi_{I})
\geq m(\det A_{II})^{1/m} f_{I}(\pi_{I}x)\geq c' g(x).
\]
All the other summands have nonnegative $L_{A}$. Thus $J=(c')^{-1}\sum_{I}J_{I}\circ\pi_{I}$ has the required bounds.
\end{proof}

For each $I$, take the mixed norm $L^{q}_{z}L^{\infty}_{y}$ using the zero extension from $B_{1/2}$. The exponents below satisfy $2-d/q=(k-1)d/q$ at $q=q_{*}$.

\begin{lemma}\label{thm3.3}
Let $u\in C^{0}(B_{1})\cap W_{\mathrm{loc}}^{2,q_{*}}(B_{1})$ be subharmonic. For each $I$, as $\varepsilon\downarrow0$,
\[
\varepsilon^{d/q}\|D^{2}U_{\varepsilon}\|_{L^{q}_{z}L^{\infty}_{y}}=o(1),\quad
\|U_{\varepsilon}-u\|_{L^{q}_{z}L^{\infty}_{y}}=o(\varepsilon^{2-d/q}).
\]
\end{lemma}

\begin{proof}
Extend $H=\|D^{2}u\|$ from $B_{3/4}$ by zero, and denote the average of $H$ over $B_{t}(x)$ by $M_{t}H(x)$. For $t<1/8$, integration in $y$ and H\"older's inequality give
\[
M_{t}H(z,y)\leq C t^{-d/q-m}\int_{B_{t}^{m}(z)}
\|H(z',\cdot)\|_{L^{q}(\mathbb{R}^{d})}\,dz'.
\]
Taking the supremum over $y$ and the $L^{q}$ norm in $z$ gives $t^{d/q}\|M_{t}H\|_{L^{q}_{z}L^{\infty}_{y}}\leq C\|H\|_{L^{q}}$, also for $q=1$. The left-hand side tends to zero for bounded $H$ because the projection of $B_{1/2}$ has finite measure. Applying the bound to $H-\min\{H,N\}$ and sending $N$ to infinity proves the same limit for every $H\in L^{q}$. Since $\|D^{2}U_{\varepsilon}\|\leq C M_{\varepsilon}H$ on $B_{1/2}$, this proves the first assertion.

By subharmonicity and radial averaging, $U_{\varepsilon}(x)-u(x)$ lies between zero and the increment of the spherical mean from radius zero to $\varepsilon$. The derivative of that mean at radius $t$ is $(t/n)M_{t}(\Delta u)(x)$. The integrated identity follows by smooth approximation in $C^{0}\cap W^{2,q}$ between two positive radii and then by letting the smaller radius tend to zero. Thus, by Minkowski's inequality,
\[
\|U_{\varepsilon}-u\|_{L^{q}_{z}L^{\infty}_{y}}
\leq C\int_{0}^{\varepsilon}t\|M_{t}H\|_{L^{q}_{z}L^{\infty}_{y}}\,dt
=o(\varepsilon^{2-d/q}),
\]
where $d/q=2/k<2$.
\end{proof}

\subsection{The Sobolev endpoint}

\begin{lemma}\label{thm3.4}
Let $2\leq k\leq n$ and let $u\in C^{0}(\Omega)\cap\Phi^{k}(\Omega)$ satisfy
$\mu_{k}[u]=f\,dx$ in $\Omega$, where $f\in C^{0}(\Omega)$ and $f>0$.
Suppose that $u\in W_{\mathrm{loc}}^{2,q_{*}}(\Omega)$ with $q_{*}=k(k-1)/2$.
Then $H_{B}u(x_{0})>u(x_{0})$ for every ball $B=B_{R}(x_{0})\Subset\Omega$.
\end{lemma}

\begin{proof}
By comparison on nested balls it suffices to consider arbitrarily small balls about $x_{0}$. Rescale by $r^{-2}u(x_{0}+rx)$ to take $\overline B_{2}\subset\Omega$ and $v=H_{B_{1/2}}u$. Suppose $v(0)=u(0)$. Set $h=U_{\varepsilon}-u$, $\delta=\varepsilon^{2-2/k}$ and $\lambda=\min_{\overline B_{2}}f$. By \eqref{eqn3.4}, $h\geq0$ and $\sigma_{k}(D^{2}U_{\varepsilon})\geq\lambda$ on $\overline B_{1}$ for small $\varepsilon$.

With $d,m,q$ as above, set $p=q/m$, $R=3/8$, and choose $0\leq\chi\leq1$ in $C_{c}^{\infty}(B_{1/2})$ equal to one on $\overline B_{R}$. Put
\[
g=\chi\big[(\delta^{1/d}\|D^{2}U_{\varepsilon}\|)^{p}+(h/\delta)^{p}\big],
\quad \tau_{\varepsilon}=\Big(\varepsilon+\sum_{I}\|g_{I}\|_{L^{m}}\Big)^{1/2}.
\]
Since $pm=q$ and $\delta^{1/d}=\varepsilon^{d/q}$, Lemma \ref{thm3.3} gives $\tau_{\varepsilon}\to0$ by $\|\sup_{y}F^{p}\|_{L^{m}_{z}}=\|\sup_{y}F\|_{L^{q}_{z}}^{p}$ for $F\geq0$. Lemma \ref{thm3.2} applied to $g/\tau_{\varepsilon}$ gives a smooth convex $J\leq0$ on $B_{1}$ with
\[
\|J\|_{L^{\infty}(B_{1})}\leq C\tau_{\varepsilon},\quad
L_{A}J\geq g/\tau_{\varepsilon}\quad\text{in }B_{1/2}\quad(A\in\mathcal{K}_{k}).
\]
Convexity gives $\|\nabla J\|_{L^{\infty}(B_{R})}\leq C\tau_{\varepsilon}$ as well.

Set $s=1-|x|^{2}/R^{2}+J$ and $D=2+(U_{\varepsilon}-v)/\delta$. Write $t_{+}=\max\{t,0\}$. Since $s\leq0$ on $\partial B_{R}$, $s(0)\to1$, and $D(0)\geq2$, the maximum $N$ of $s_{+}D_{+}^{p}$ on $\overline B_{R}$ is at least one for small $\varepsilon$. It is attained at an interior point $x$ with $s,D>0$, so
\[
Q=U_{\varepsilon}+2\delta-\delta(N/s)^{1/p}
\]
touches $v$ from below at $x$. By the contact property in Section 3.1 and duality, $L_{A}Q(x)\leq0$ for some $A\in\mathcal{K}_{k}$.

At $x$ we have $(N/s)^{1/p}=D$. Differentiating $Q$ with $N$ fixed and using $s\leq1$, $|\nabla s|\leq C$ and $\operatorname{tr} A=1$ gives
\begin{align*}
0\geq L_{A}Q
&=L_{A}U_{\varepsilon}+\frac{\delta D}{ps}
\left(L_{A}s-(1+1/p)\frac{(\nabla s)^{t} A\nabla s}{s}\right)\\
&\geq L_{A}U_{\varepsilon}+\frac{\delta D}{ps}
\left(\frac{g}{\tau_{\varepsilon}}-\frac{C}{s}\right).
\end{align*}
Since $N=sD^{p}\geq1$, $v\geq u$ and $\chi(x)=1$, we have $s^{-1}\leq D^{p}\leq(2+h/\delta)^{p}\leq C(1+g)$. Also $L_{A}U_{\varepsilon}\geq0$ by admissibility. Hence $g/\tau_{\varepsilon}\leq C(1+g)$, which implies $g(x)\leq2C\tau_{\varepsilon}$ for small $\varepsilon$. Both $D$ and $s^{-1}$ are now bounded at $x$, and the last differential inequality gives $L_{A}U_{\varepsilon}(x)\leq C_{0}\delta$.

But $g(x)\to0$ also gives $\delta\|H\|^{k-1}\to0$ for $H=D^{2}U_{\varepsilon}(x)$. Hence
\[
\sup_{0\leq t\leq1}\big|\sigma_{k}(H-t(C_{0}+1)\delta I)-\sigma_{k}(H)\big|
\leq C\delta(\|H\|+\delta)^{k-1}=o(1).
\]
Since $\sigma_{k}(H)\geq\lambda>0$, the segment stays in $\Gamma_{k}$, so $L_{A}U_{\varepsilon}(x)\geq(C_{0}+1)\delta$. This contradiction proves $v(0)>u(0)$.
\end{proof}

For $k=2$, Theorem \ref{thm1.1} supplies the endpoint hypothesis. The positive constant-density comparison solution used by the author \cite[Section 3.1, (3.1)--(3.3)]{LyuSigma2Interior} also gives separation, since it lies below $H_{B}u$.

\section{Interior $W^{2,p}$ regularity}

Throughout this section $2\leq k\leq n$. For a ball $B\Subset\Omega$, write $H_{B}u$ for the solution of \eqref{eqn3.2} and $U_j$ for the smooth Dirichlet approximations from Section 3. A comparison function $W$ from \eqref{eqn3.3} provides a weight for estimating $\Delta U_j$ where $H_Bu>u$.

\begin{theorem}\label{thm4.1}
Let $u\in C^{0}(\Omega)\cap\Phi^{k}(\Omega)$ satisfy \eqref{eqn1.4} and let $B\Subset\Omega$ be a ball. For $1<p<\infty$ and any nonempty open set $D\Subset B\cap\{H_{B}u>u\}$,
\[
\|D^{2}U_{j}\|_{L^{p}(D)}\leq C
\quad\text{for all sufficiently large }j.
\]
In particular, $u\in W_{\mathrm{loc}}^{2,p}(B\cap\{H_{B}u>u\})$.
\end{theorem}

Choose $W$ as in \eqref{eqn3.3} with $K=\overline D$ and a closed boundary collar of $B$ on which $W<u$. The constant in Theorem \ref{thm4.1} depends only on $n,k,p$, the radius of $B$, $\|u\|_{L^{\infty}(B)}$, positive lower and upper bounds and a modulus of continuity for $f$ on $\overline B$, $\|W\|_{C^{2}(\overline B)}$, $\min_{\overline D}(W-u)$ and the width of the collar.

\subsection{Stability under density perturbations}

We bound the logarithm of the ratio of the Laplacians and a weighted difference of the solutions by the same density error.

\begin{lemma}\label{thm4.2}
Let $U,v\in C^{2}(\overline{B})\cap C^{\infty}(B)$ be admissible on a ball $B$ of radius $R$, with $U=v$ on $\partial B$.
If $g=G(D^{2}U)$ and $h=G(D^{2}v)$ satisfy $0<m\leq g,h\leq M$, then for every $a\geq2(k-1)$,
\begin{equation}\label{eqn4.1}
\left\|\log\frac{\Delta v}{\Delta U}\right\|_{L^{2}(B)}
+R^{-2}\big\|(\Delta U)^{1/a}(v-U)\big\|_{L^{2}(B)}
\leq C(n,k,m,M)\,\|h-g\|_{L^{2}(B)} .
\end{equation}
\end{lemma}
\begin{proof}
Set $A=D^{2}U$, $B'=D^{2}v$, $X=A/g$, $Y=B'/h$, and $t=\operatorname{tr} Y/\operatorname{tr} X$.
The positive degree-one function $Q=(\sigma_{k}/\operatorname{tr})^{1/(k-1)}$ is concave, so $DQ(X)[Y]\geq Q(Y)$.
Logarithmic differentiation of $\sigma_{k}=(\operatorname{tr})Q^{k-1}$ and $\sigma_{k}(X)=\sigma_{k}(Y)=1$ give
\begin{align*}
D\sigma_{k}(X)[Y]-k
&=t+(k-1)\frac{DQ(X)[Y]}{Q(X)}-k\\
&\geq t+(k-1)t^{-1/(k-1)}-k\\
&=:F(t)\geq0.
\end{align*}
Here $F(t)\geq0$ by the arithmetic-geometric mean inequality. Interchanging $X$ and $Y$ gives $F(t)+F(t^{-1})\geq t+t^{-1}-2\geq(\log t)^{2}$.
For the integral estimate we use the symmetric expression of Chen-Zhou-Zhu for $k=2$ \cite[Proposition 3.8]{CZZ} and its general-$k$ form given by the author \cite[Lemma 2.3]{LyuKGlobal} (see also the author \cite[Section 2.2]{LyuSigma2Global}). Set
\[
E=D\sigma_{k}(A)[B']-kg^{k-1}h+D\sigma_{k}(B')[A]-kh^{k-1}g .
\]
By homogeneity,
\[
E=g^{k-1}h\big(D\sigma_{k}(X)[Y]-k\big)
+h^{k-1}g\big(D\sigma_{k}(Y)[X]-k\big).
\]
Since both differences are nonnegative, $E\geq m^{k}(F(t)+F(t^{-1}))\geq m^{k}(\log t)^{2}$.

Put $w=v-U$. Euler's identity and integration by parts with the divergence-free Newton tensors give
\[
\int_{B} E=k\int_{B}(g-h)(g^{k-1}-h^{k-1})-
\int_{\partial B}w_{\nu}\big(\sigma_{k-1}(B'|_{T\partial B})
-\sigma_{k-1}(A|_{T\partial B})\big)\, dS .
\]
This identity follows by integrating first on concentric interior balls and passing to the boundary. Since $w=0$ on $\partial B$, its tangential Hessian is $(w_{\nu}/R)I$, where $\nu$ is outward. Both tangential blocks lie in $\Gamma_{k-1}$, so ellipticity makes the boundary integrand nonnegative. Thus $\int_{B}E\leq C\|g-h\|_{L^{2}(B)}^{2}$.
Use $\log(\Delta v/\Delta U)=\log t+\log(h/g)$ and the Lipschitz bound for $\log$ on $[m,M]$ to obtain the first term in \eqref{eqn4.1}.

For the second term, center $B$ at zero and put $T=\Delta U$. Let $\mathcal A$ be the average of the positive definite, divergence-free tensors $D\sigma_k((1-t)D^{2}U+tD^{2}v)$ over $0\leq t\leq1$. Then
\[
\partial_{i}(\mathcal{A}^{ij}w_{j})=h^{k}-g^{k},\qquad w|_{\partial B}=0.
\]
By concavity $G((1-t)D^{2}U+tD^{2}v)\geq m$, and for $t\leq1/2$ the trace of the matrix is at least $T/2$.
The Newton-Maclaurin inequality gives
\[
\sigma_{k-1}(H)^{k-1}
\geq c(n,k)\,\sigma_{1}(H)\sigma_{k}(H)^{k-2}.
\]
Since $\operatorname{tr}D\sigma_k(H)=(n-k+1)\sigma_{k-1}(H)$, integration on $0\leq t\leq1/2$ and the Newton-Maclaurin inequality yield $\operatorname{tr}\mathcal{A}\geq c(1+T^{1/(k-1)})$.

Put $I=\int_{B} w^{2}\operatorname{tr}\mathcal{A}$ and $\mathcal{E}=\int_{B}\mathcal{A}^{ij}w_{i}w_{j}$.
Since $\partial_{i}(\mathcal{A}^{ij}x_{j})=\operatorname{tr}\mathcal{A}$, integration by parts and the Cauchy-Schwarz inequality give
\[
I=-2\int_{B} w\mathcal{A}^{ij}x_{j}w_{i}
\leq2R I^{1/2}\mathcal{E}^{1/2},
\]
\[
\mathcal{E}=-\int_{B} w(h^{k}-g^{k})
\leq C I^{1/2}\|h-g\|_{L^{2}(B)}.
\]
Thus $I^{1/2}\leq CR^{2}\|h-g\|_{L^{2}(B)}$.
The second estimate in \eqref{eqn4.1} follows since $T^{2/a}\leq1+T^{1/(k-1)}\leq C\operatorname{tr}\mathcal A$.
\end{proof}

\subsection{A weighted Hessian bound}

Set $a:=\max\{4,2(k-1)\}$. The bound $a\geq4$ controls the third derivatives in the proof below, and $a\geq2(k-1)$ permits the use of \eqref{eqn4.1}.

\begin{lemma}\label{thm4.3}
Let $V\Subset V'\Subset B$. Suppose $v\in C^{\infty}(B)$ is admissible, $\|v\|_{L^{\infty}(B)}\leq A$ and $0<m\leq G(D^{2}v)\leq M$. Let $W\in C^{\infty}(\overline B)$ be $k$-convex with $\{W>v\}\Subset V$ and set $Y=(W-v)_{+}^{a}\Delta v$.
Given $b\in C^{\infty}(V')$, let $P$ satisfy
\[
P\geq C\bigl(1+\|\nabla b\|_{L^{\infty}(V')}
+\|D^{2}b\|_{L^{\infty}(V')}^{1/2}\bigr),
\]
where $C>0$ depends only on $n,k,m,M,A$, $\operatorname{dist}(V',\partial B)$ and $\|W\|_{C^{2}(V')}$.
If $G(D^{2}v)=b$ on $\{Y>2P\}$, then
\begin{equation}\label{eqn4.2}
\sup_{B}Y\leq C_{0}P,
\end{equation}
where $C_{0}>0$ depends only on $n$ and $k$.
\end{lemma}
\begin{proof}
Put $\rho=W-v$, $O=\{\rho>0\}$, $L=\|\nabla b\|_{L^{\infty}(V')}$ and $J=\|D^{2}b\|_{L^{\infty}(V')}$.
With $\delta=\operatorname{dist}(V',\partial B)$, the interior gradient estimate of Liu-Ma-Zang \cite[Corollary 1.4]{LMZ} on balls of radius $\delta/2$ gives
\[
\sup_{V'}|\nabla v|\leq C(n,k)(A/\delta+\delta M)
\]
for $k<n$ (see also He \cite[Theorem 1.1 and Proposition 2.1]{He}). For $k=n$, convexity gives $\sup_{V'}|\nabla v|\leq2A/\delta$.
Thus $K:=1+\sup_{V'}(|\nabla v|^{2}+|\nabla W|^{2})$ is bounded by the stated data.

If $O$ is empty, $Y=0$. Otherwise put $\theta=(8K)^{-1}$ and maximize
\[
Q(x,e)=\rho^{a} e^{\theta|\nabla v|^{2}/2} v_{ee}
\]
over $\overline O$ and unit vectors $e$, as in Pogorelov \cite{Pogorelov} and  Chou-Wang \cite[Theorem 4.1]{CW}.
The maximum is positive and lies in $O$. Since $\lambda_{\max}(D^{2}v)\leq\Delta v\leq n\lambda_{\max}(D^{2}v)$, we have $\sup_B Y\leq n\max Q$. If $Y\leq2P$ at the maximum, then $\max Q\leq2e^{1/16}P$. Otherwise we may differentiate $G(D^{2}v)=b$ there (see Chen-Zhou-Zhu \cite[Section 3.3]{CZZ}; see also the author \cite[Lemma 2.4]{LyuSigma2Interior}).

Choose an eigenframe at the maximum with $e=e_{1}$ and $\lambda_{1}\geq\cdots\geq\lambda_{n}$. Put $F_{i}=G^{ii}(D^{2}v)$ and $A_{1}=D^{2}v_{1}$. The first derivative condition gives
\[
\frac{v_{11i}}{\lambda_{1}}=-a\frac{\rho_{i}}{\rho}-\theta v_{i}\lambda_{i}.
\]
By concavity and homogeneity, $DG(D^{2}v)[D^{2}W]\geq0$.
Applying the linearized operator to $\log Q(\cdot,e_{1})$ and differentiating $G(D^{2}v)=b$ in the $e_{1}$ direction twice, we get
\begin{equation}\label{eqn4.3}
0\geq\theta\sum_i F_i\lambda_i^2-a\sum_i F_i\frac{\rho_i^2}{\rho^2}
-\frac{aM}{\rho}-\theta K^{1/2} L-\frac{J}{\lambda_1}-\frac{D^{2}G(D^{2}v)[A_1,A_1]}{\lambda_1}
-\sum_i F_i\frac{v_{11i}^2}{\lambda_1^2}.
\end{equation}
Fix $\varepsilon>0$ depending only on $n,k$ with $(n-k)\varepsilon\leq1/3$.
If $\lambda_k\geq\varepsilon\lambda_1$, then $\sum_iF_i\lambda_i^2\geq F_k\lambda_k^2\geq c\lambda_1^2\sum_iF_i$ by Chou-Wang \cite{CW}.
The first derivative identity gives
\[
\sum_iF_i\frac{v_{11i}^2}{\lambda_1^2}
\leq 2a^2\sum_i F_i\frac{\rho_i^2}{\rho^2}
+2\theta^2K\sum_iF_i\lambda_i^2.
\]
The last term is absorbed by the first term of \eqref{eqn4.3}, and $-D^{2}G[A_1,A_1]\geq0$.

If $\lambda_k<\varepsilon\lambda_1$, we use
\[
\lambda_1\sigma_{k-2}(\lambda|1i)\geq\tfrac34\sigma_{k-1}(\lambda|i),
\quad \forall i\geq2.
\]
Here $\lambda|I$ denotes deletion of the entries indexed by $I$. Write $\mu=\lambda|1i$ in decreasing order. By hyperplane restriction, $\mu\in\Gamma_{k-2}$, and
$\sigma_{k-1}(\lambda|i)=\lambda_1\sigma_{k-2}(\mu)+\sigma_{k-1}(\mu)$.
If $\sigma_{k-1}(\mu)\leq0$ the inequality follows immediately. Otherwise $\mu\in\Gamma_{k-1}$ and $0<\mu_{k-1}\leq\lambda_k$.
The quotient $\sigma_{k-1}/\sigma_{k-2}$ increases in each entry on $\Gamma_{k-1}$ (see Trudinger \cite{TrudQuotient}). Replacing the last $n-k$ entries of $\mu$ by $\mu_{k-1}$, and then letting the first $k-2$ entries tend together to $+\infty$, gives
\[
\frac{\sigma_{k-1}(\mu)}{\sigma_{k-2}(\mu)}
\leq(n-k)\mu_{k-1}\leq\lambda_1/3 .
\]
For $k=2$ the first group of entries is empty, and for $k=n$ one always has $\sigma_{k-1}(\mu)=0$. This proves the claimed inequality.
The diagonal part of $D^2G$ is nonpositive. Keeping the off-diagonal terms with indices $1,i$ as in Chou-Wang \cite{CW} therefore gives
\[
-\frac{D^2G(D^2v)[A_1,A_1]}{\lambda_1}
\geq \frac32\sum_{i\geq2}F_i\frac{v_{11i}^2}{\lambda_1^2}.
\]
For $i\geq2$ the first derivative identity gives
\[
a\frac{\rho_i^2}{\rho^2}
\leq\frac{2}{a}\left(\frac{v_{11i}^2}{\lambda_1^2}
+\theta^2v_i^2\lambda_i^2\right).
\]
Since $2/a\leq1/2$, the third derivatives for $i\geq2$ cancel in \eqref{eqn4.3}. For $i=1$ use the first-case bound, and absorb the terms with factor $\theta^2K$ as before.

In the first case put $D_0=\sum_i F_i$, and in the second put $D_0=F_1$. Since $|\nabla\rho|^2\leq2K$, both cases yield
\[
c\theta D_0\lambda_1^2\leq CKD_0/\rho^2+C/\rho
+C\theta K^{1/2} L+CJ/\lambda_1 .
\]
If $\rho\lambda_1\leq CK$, the boundedness of $\rho$ gives $\max Q\leq C$.
Otherwise the first term on the right is absorbed into the left. By \eqref{eqn2.10}, $F_1\lambda_1\geq m/(nk)$, so
\[
\lambda_1\leq C(K/\rho+K^{1/2}L+\sqrt{KJ}).
\]
Multiplication by $\rho^a e^{\theta|\nabla v|^2/2}$ gives $\max Q\leq C(1+L+\sqrt J)$. Taking $P$ above this bound and combining the two alternatives proves \eqref{eqn4.2}.
\end{proof}

\subsection{Level-set decay}

\begin{proof}[Proof of Theorem \ref{thm4.1}]
Center $B$ at the origin and let $R$ be its radius. Write $g=f^{1/k}$ and $g_{j}=G(D^{2}U_{j})$. Choose $m,M$ from the density bounds so that $0<2m\leq g\leq M/2$ on $\overline B$. The modulus of $g$ is controlled by that of $f$ and the lower density bound.

Set $\gamma=\min_{\overline D}(W-u)>0$. For large $j$, uniform convergence gives $W-U_j\geq\gamma/2$ on $\overline D$, $U_j>W$ on the collar and $\|U_j\|_{L^\infty(B)}\leq1+\|u\|_{L^\infty(B)}$. In particular, $\{W>U_j\}\Subset V$ for concentric balls $V\Subset V'\Subset B$ chosen from the collar width.

Fix a smooth nondecreasing cutoff $\eta$ that is zero on $(-\infty,1]$ and one on $[2,\infty)$.
Let $0<d\leq\min\{1,m,d_{*}\}$, where $d_{*}$ is the threshold in Proposition \ref{thmA.1}. Extend $g$ to $\mathbb R^{n}$ by nearest-point projection onto $\overline B$. Mollify the extension at a scale $\varepsilon>0$ such that the error on $\overline B$ is at most $d/8$, and subtract $d/2$ to obtain $b\in C^{\infty}(\overline B)$.
For large $j$ with $\|g_{j}-g\|_{L^\infty(B)}\leq d/8$ we have
\[
m\leq b\leq g_{j}\leq M,\quad
\|g_{j}-b\|_{L^{\infty}(B)}\leq d.
\]
The projection is $1$-Lipschitz, so a common modulus of $g$ permits a common choice of $\varepsilon$. The first and second derivatives of $b$ are bounded in terms of $n,M,\varepsilon$ and the fixed mollifier.

Fix such a $j$ and write $U=U_{j}$, $T=\Delta U$, $\rho=(W-U)_{+}$ and $X=\rho^{a}T$. For each $P>0$, set
\[
\Theta=\frac{P}{2T}+\rho^{a}\eta(4X/P).
\]
The cutoff makes the second term zero near $\{\rho=0\}$, so $\Theta$ is smooth and positive. The first term is at least $\rho^{a}$ for $X\leq P/2$, and the cutoff equals one for $X>P/2$. Hence $\Theta\geq\rho^{a}$, while $\Theta T/P=1/2$ on $\{X\leq P/4\}$.
Proposition \ref{thmA.1} gives an admissible $v\in C^{2}(\overline B)\cap C^{\infty}(B)$ with
\[
h=G(D^{2}v)=g_{j}\exp\!\left[-\log(g_{j}/b)\eta(\Theta\Delta v/P)\right]
\quad\text{in }B,\qquad v=U\quad\text{on }\partial B .
\]
For a density cutoff in a related estimate, see the author \cite[Appendix A, (A.9)]{LyuSigma2Higher}.
Comparison and subharmonicity give $U\leq v\leq\max_{\partial B}U$ and $\|v\|_{L^{\infty}(B)}\leq1+\|u\|_{L^{\infty}(B)}$.

Write $S=\Delta v$ and $Y=(W-v)_{+}^{a}S$. Then $\{W>v\}\Subset V$ and $Y\leq\Theta S$. On $\{Y>2P\}$ the cutoff is one, so $h=b$ there. Lemma \ref{thm4.3} gives $Y\leq C_{0}P$ whenever $P\geq P_{b}$, where $C_0\geq1$ depends only on $n,k$ and
\[
P_{b}=\max\!\left\{1,R^{2a},
C\big(1+\|\nabla b\|_{L^\infty(V')}+\|D^{2}b\|_{L^\infty(V')}^{1/2}\big)\right\}
\]
with the constant $C$ from Lemma \ref{thm4.3}.
This choice gives $R^{4}P^{-2/a}\leq1$ for $P\geq P_b$, uniformly in $j$.

Let $N(t)=|\{x\in B:X(x)>t\}|$ and denote
\[
E=\{h\ne g_{j}\},\quad F=\{S>2T\}\cup\{T>2S\}.
\]
On $F$ we have $|\log(S/T)|\geq\log 2$. Lemma \ref{thm4.2} and the support of the density change give
\[
|F|\leq C\|h-g_{j}\|_{L^{2}(B)}^{2}\leq Cd^{2}|E|.
\]
On $\{X\leq P/4\}\setminus F$, we have $\Theta S/P\leq1$, and hence $h=g_{j}$. Thus $E\subset\{X>P/4\}\cup F$ and
\[
|E|\leq 2N(P/4)\quad\text{if }Cd^{2}\leq\tfrac12 .
\]

On $B\setminus F$ we have $T\leq2S$. Since $\rho\leq (W-v)_{+}+|v-U|$,
\[
X^{1/a}\leq(2Y)^{1/a}+|v-U|T^{1/a}.
\]
Thus, for $P\geq P_b$, the set $\{X>4C_{0}P\}\setminus F$ is contained in $\{|v-U|T^{1/a}>cP^{1/a}\}$ with $c=(4C_{0})^{1/a}-(2C_{0})^{1/a}>0$. By \eqref{eqn4.1} and Chebyshev's inequality, its measure is at most $CR^{4}P^{-2/a}\|h-g_{j}\|_{L^{2}(B)}^{2}$. Adding the bound for $|F|$ gives
\begin{equation}\label{eqn4.4}
N(4C_{0}P)\leq C_{*}d^{2}N(P/4)\quad(P\geq P_{b}).
\end{equation}
The constant $C_{*}$ depends only on $n,k,m,M$.
For related distribution estimates see the author \cite[Section 3.3]{LyuSigma2Interior} and, with boundary weights, the author \cite[Section 2.4]{LyuSigma2Global,LyuKGlobal}.

Choose $d$ so that, in addition to the preceding smallness conditions,
$C_{*}d^{2}(16C_{0})^{p}\leq\tfrac12$. This choice uses only $n,k,p,m,M$ and the fixed cutoff. With $b$ and $P_b$ fixed, split $p\int_{0}^{\infty}t^{p-1}N(t)\,dt$ at $t=4C_{0}P_b$. In the upper integral apply \eqref{eqn4.4} with $P=t/(4C_{0})$ and change variables to $t/(16C_{0})$. This gives
\[
\int_{B}X^{p}
\leq |B|(4C_{0}P_b)^{p}+C_{*}d^{2}(16C_{0})^{p}\int_{B}X^{p}.
\]
For each $j$, $X$ is bounded, so absorption gives $\|X\|_{L^{p}(B)}\leq (2|B|)^{1/p}4C_{0}P_b$.

On $D$, $\rho\geq\gamma/2$ and $\|D^{2}U\|\leq T$ by $k$-admissibility. Hence
\[
\|D^{2}U_{j}\|_{L^{p}(D)}
\leq \left(\frac{2}{\gamma}\right)^{a}(2|B|)^{1/p}4C_{0}P_b .
\]
The boundary gap and convergence of $g_j,U_j$ determine how large $j$ must be. The bound itself has the stated dependence. Weak compactness and uniform convergence give $D^{2}u\in L^{p}(D)$. Interior interpolation yields $u\in W^{2,p}_{\mathrm{loc}}(D)$, and $D$ was arbitrary.
\end{proof}

\subsection{The interior estimate}

Together with Lemmas \ref{thm3.1} and \ref{thm3.4}, Theorem \ref{thm4.1} gives the following local characterization of strict separation.

\begin{proposition}\label{thm4.4}
Let $2\leq k\leq n$ and let $u\in C^{0}(\Omega)\cap\Phi^{k}(\Omega)$ satisfy
\[
\mu_{k}[u]=f\,dx \quad\text{in }\Omega,
\qquad f\in C^{0}(\Omega),\quad f>0.
\]
Fix $x_{0}\in\Omega$. Then the following conditions are equivalent.
\begin{enumerate}
\item[(i)] $H_{B}u(x_{0})>u(x_{0})$ for some ball $B\Subset\Omega$ centered at $x_{0}$.
\item[(ii)] $u\in W^{2,q_{*}}(B_{r}(x_{0}))$ for some $r>0$, where $q_{*}=k(k-1)/2$.
\item[(iii)] For some $r>0$, $u\in W^{2,p}(B_{r}(x_{0}))$ for every $1<p<\infty$.
\item[(iv)] $u\in C^{1,\beta}(B_{r}(x_{0}))$ for some $r>0$ and some $1-2/k<\beta<1$.
\end{enumerate}
\end{proposition}
\begin{proof}
Under (i), choose $r>0$ with $\overline B_r(x_0)\subset B\cap\{H_Bu>u\}$. The local conclusion of Theorem \ref{thm4.1} gives (iii). In (iii), taking $p=\max\{2,q_*\}$ gives (ii) by H\"older's inequality, and taking $p>nk/2$ gives (iv) by Morrey's inequality. Lemma \ref{thm3.4} gives (i) from (ii), and Lemma \ref{thm3.1} gives (i) from (iv).
\end{proof}

Proposition \ref{thm4.4} gives local $W^{2,p}$ regularity for each solution. To obtain the bound in Theorem \ref{thm1.2}, we must control the dependence of the estimate on the comparison functions. We choose these functions for a limit solution. Their strict gaps persist along the convergent subsequence, so Theorem \ref{thm4.1} applies with uniform constants.

\begin{lemma}\label{thm4.5}
Let $2\leq k\leq n$ and let $B_r\Subset B_R\Subset\Omega$ be concentric balls. Suppose that $u_i\in C^0(\Omega)\cap\Phi^k(\Omega)$ satisfy $\mu_k[u_i]=f_i\,dx$ in $\Omega$. Assume that the $u_i$ are uniformly bounded on $B_R$ and that the $f_i$ have common positive lower and upper bounds and a common modulus of continuity on $B_R$. For $k=2$, these hypotheses suffice. For $k\geq3$, assume in addition that either
\[
\sup_i\|u_i\|_{W^{2,q_*}(B_R)}<\infty
\quad\text{or}\quad
\sup_i[\nabla u_i]_{C^{0,\beta}(B_R)}<\infty
\quad\text{for some fixed }\beta_*<\beta<1.
\]

After passing to a subsequence, there are finitely many pairs of concentric balls $D_\ell\Subset B_\ell\Subset B_R$ such that each $D_\ell$ has radius at most one and the half balls $\tfrac12 D_\ell$ cover $\overline B_r$. For every $1<p<\infty$, the smooth Dirichlet approximations $U_{i,j,\ell}$ of $u_i$ on $B_\ell$ satisfy
\[
\|D^2U_{i,j,\ell}\|_{L^p(D_\ell)}\leq C_{p,\ell}
\quad\text{for all }i\text{ and }j\geq j_0(i,p,\ell).
\]
\end{lemma}

\begin{proof}
Fix $r<s<R$. The interior H\"older estimate of Trudinger \cite{TrudWeak}, in the form given by Trudinger-Wang \cite[(6.4)]{TW2}, gives equicontinuity of the $u_i$ on $\overline B_s$. The $f_i$ also have a common modulus of continuity. Passing to a subsequence, we have $u_i\to u$ and $f_i\to f$ uniformly on $\overline B_s$. The limit $u$ is $k$-convex, and weak continuity of Hessian measures gives $\mu_k[u]=f\,dx$ in $B_s$, with $f$ continuous and strictly positive.

The limit satisfies one of the separation criteria of Section 3. For $k\geq3$ in the Sobolev case, weak compactness gives $u\in W^{2,q_*}(B_s)$ since $q_*>1$. In the H\"older case, difference quotients of fixed length less than $R-s$, together with the common gradient seminorm bound, give a uniform bound for $\nabla u_i$ on $\overline B_s$. A further subsequence of the gradients converges uniformly to $\nabla u$, which inherits the same H\"older bound. For $k=2$, Theorem \ref{thm1.1} applied to the limit equation gives $u\in W^{2,1}_{\mathrm{loc}}(B_s)$, and $q_*=1$.

For each $x\in\overline B_r$, take a ball $B\Subset B_s$ centered at $x$. Lemma \ref{thm3.1} or Lemma \ref{thm3.4} gives $H_Bu(x)>u(x)$. By continuity and \eqref{eqn3.3}, there are a smaller concentric ball $D$ of radius at most one and a smooth admissible $W$ on $\overline B$ such that $W>u$ on $\overline D$ and $W<u$ on $\partial B$. Choose finitely many of the half balls $\tfrac12 D$ to cover $\overline B_r$, and denote the corresponding balls and functions by $D_\ell$, $B_\ell$ and $W_\ell$.

For each $\ell$, the inequality $W_\ell<u$ persists on a closed boundary collar of $B_\ell$. Since the cover is finite, uniform convergence allows us to discard finitely many indices so that $W_\ell-u_i$ on $\overline D_\ell$ and $u_i-W_\ell$ on the collar have a common positive lower bound for all $i$ and $\ell$. By comparison, $H_{B_\ell}u_i\geq W_\ell>u_i$ on $\overline D_\ell$. On each ball the comparison function, the gap bound and the collar are now fixed, so Theorem \ref{thm4.1} gives the asserted estimate uniformly in $i$.
\end{proof}

\begin{proof}[Proof of Theorem \ref{thm1.2}]
Fix $1<p<\infty$. If the asserted bound failed, there would be a sequence $u_i$ with common bounds for the stated data and $\|u_i\|_{W^{2,p}(B_r)}\to\infty$. Take the subsequence and finite cover given by Lemma \ref{thm4.5}. For fixed $i$ and $\ell$, the uniform convergence $U_{i,j,\ell}\to u_i$ as $j\to\infty$ identifies the distributional limit of the Hessians with $D^2u_i$. Weak compactness in $L^p(D_\ell)$ and lower semicontinuity therefore give
\[
\|D^2u_i\|_{L^p(D_\ell)}\leq C_{p,\ell}.
\]
Interior interpolation and the common supremum bound give a uniform $W^{2,p}$ bound on each $\tfrac12 D_\ell$. These finitely many half balls cover $\overline B_r$, contradicting the divergence of the norms.
\end{proof}

\section{Interior $W^{2,\infty}$ and $C^{2,\alpha}$ regularity}

\subsection{Square-root Dini densities}

We combine a lower bound for the average of $\Delta u$ with the $L^p$ estimate from Theorem \ref{thm4.1} to obtain an interior Hessian bound.

\begin{lemma}\label{thm5.1}
Let $2\leq k\leq n$ and $B=B_R(x_0)$ with $0<R\leq1$. Let $u$ be smooth and $k$-admissible near $\overline B$, and suppose that
\[
G(D^2u)=g\quad\text{in }B,
\]
where $0<m\leq g\leq M_0$ and $g$ has a nondecreasing modulus of continuity $\omega$ on $B$ satisfying $\Psi_\omega(1)<\infty$, with $\Psi_\omega$ as in \eqref{eqn1.5}.

Set $M=\Delta u(x_0)$. If $M\geq1$ and $\sup_B\Delta u\leq AM$ for some $A\geq1$, then there exists $y\in B_{R/4}(x_0)$ such that
\begin{equation}\label{eqn5.1}
\frac{1}{|B_h|}\int_{B_h(y)}\Delta u\,dx\geq c_1M,
\quad h=c_0RM^{-k/2},
\end{equation}
where $0<c_0<1/4$ and $c_1>0$ depend only on $n,k,m,M_0,A$ and $\omega$.
\end{lemma}

\begin{proof}
Take $x_{0}=0$, set $S=\Delta u$, and truncate $\omega$ at $M_{0}$.
Let $M_{r}$ and $A_{r}$ denote spherical and ball averages with variable center. As in Chen-Zhou-Zhu \cite[Section 2.1]{CZZ}, set
\[
J_{r}w=\int_{0}^{1}M_{tr}w\,dt,\quad
V_{r}=\frac{6n}{r^{2}}(J_{r}u-u).
\]
For related estimates using the resolvent, see the author \cite[Lemma 3.5]{LyuSigma2Global} and \cite[Lemma 3.1]{LyuKGlobal}.
Using $\partial_{s} M_{s}u=(s/n)A_{s}S$, we obtain
\[
J_{r}u-u=\frac1n\int_{0}^{r}(1-s/r)sA_{s}S\,ds.
\]
Hence
\[
V_{r}=6\int_{0}^{1}t(1-t)A_{tr}S\,dt,\quad
r\partial_{r}V_{r}=6\int_{0}^{1}t(3t-2)A_{tr}S\,dt.
\]
The weights have integrals one and zero, respectively. For $w\geq0$, polar integration gives
\[
\int_{0}^{1}A_{tr}w\,dt
=\frac{n}{n-1}\int_{0}^{1}(1-t^{n-1})M_{tr}w\,dt
\leq\frac{n}{n-1}J_{r}w .
\]
Thus $0<V_{r}\leq AM$, $V_{r}\leq CJ_{r}S$ and $|r\partial_{r}V_{r}|\leq CJ_{r}|S-a|$ for every constant $a$.

For $0<r<R/4$, let
\[
N(r)=\max_{\overline B_{R/2}}\big(V_{r}(x)-32AM|x|^{2}/R^{2}\big).
\]
The bounds $V_r(0)>0$ and $V_r\leq AM$ place every maximizer $x_{r}$ in $B_{R/4}$. At $x_{r}$, put $\bar H=J_{r}D^{2}u$ and $\bar S=J_{r}S$. Then
\[
\bar H=D^{2}u(x_{r})+\frac{r^{2}}{6n}D^{2}V_{r}(x_{r})
\leq D^{2}u(x_{r})+CM(r/R)^{2}I .
\]
By admissibility, $\bar H\in\Gamma_{k}$ and $\|D^{2}u\|\leq AM$. Ellipticity and the polynomial bound for $D\sigma_k$ give
\[
\sigma_{k}(\bar H)\leq g(x_{r})^{k}+CM^{k}(r/R)^{2}.
\]

Whenever $N(r)\geq M/2$, the bound $V_{r}\leq CJ_{r}S$ gives $cM\leq\bar S\leq AM$.
Set $T=S/\bar S$, $b=1/(k-1)$ and $g_{0}=g(x_{r})$.
Concavity of $(\sigma_{k}/\sigma_{1})^{b}$ gives
$J_{r}(T^{-b}g^{1+b})\leq\sigma_{k}(\bar H)^{b}$.
On $B_{r}(x_{r})$, $(g_{0}/g)^{1+b}\leq1+C\omega(r)$. Since $0<b\leq1$ and $kb=1+b$, we obtain
\[
J_{r}(T^{-b})\leq (1+C\omega(r))
\left(\frac{\sigma_{k}(\bar H)}{g_{0}^{k}}\right)^{b}
\leq 1+C\big(\omega(r)+M^{k}(r/R)^{2}\big).
\]
Now $J_{r}T=1$ and $0<T\leq C$. On this interval $t^{-b}\geq1-b(t-1)+c(t-1)^{2}$, so
\[
J_{r}|S-\bar S|^{2}\leq CM^{2}\big(\omega(r)+M^{k}(r/R)^{2}\big).
\]
The function $N$ is locally Lipschitz and $N'(r)=\partial_{r}V_{r}(x_{r})$ at each point of differentiability, by testing the difference quotient at $x_{r}$. The derivative estimate for $V_{r}$ with $a=\bar S$ and the Cauchy-Schwarz inequality give
\[
|rN'(r)|\leq CM\big(\omega(r)^{1/2}+M^{k/2}r/R\big)
\quad\text{whenever }N(r)\geq M/2 .
\]
Since $V_r\to S$ uniformly on $\overline{B}_{R/2}$, $\lim_{r\downarrow0}N(r)\geq M$. Integrating up to the first radius with $N(r)=M/2$ gives
\[
N(r)\geq M-CM\big(\Psi_{\omega}(r)+M^{k/2}r/R\big).
\]
Choose $0<c_{0}<1/4$ so that $C(\Psi_{\omega}(c_{0})+c_{0})\leq1/4$.
The right-hand side is at least $3M/4$ for $r\leq h=c_{0}RM^{-k/2}$, so $V_{h}(x_{h})\geq N(h)\geq3M/4$.
In $V_{h}(x_{h})$, the contribution from $0<t<t_{0}$ is at most $3AMt_{0}^{2}\leq M/4$ for small $t_{0}=t_{0}(A)>0$. On $t_{0}\leq t\leq1$, $A_{th}S(x_{h})\leq t_{0}^{-n}A_{h}S(x_{h})$ by positivity. This gives \eqref{eqn5.1}.
\end{proof}

Let $D$ be a ball of radius at most one and $v$ a smooth admissible solution of $G(D^{2}v)=g$ near $\overline D$, with $0<m\leq g\leq M_0$ and square-root Dini modulus $\omega$ on $D$.
As in the author \cite[Section 4]{LyuSigma2Interior}, choose $x$ maximizing $d\,\Delta v$ on $\overline D$, where $d$ is distance to $\partial D$. Put $M=\Delta v(x)$, $R=d(x)/2$ and $p=n+nk/2$. Since $\Delta v\leq2M$ on $B_{R}(x)$, Lemma \ref{thm5.1} and Jensen's inequality give, for $M\geq1$,
\[
\int_{D}(\Delta v)^{p}\,dx
\geq cR^{n}M^{p-nk/2}
=c(RM)^{n}.
\]
For $M<1$, $d(x)M\leq1$. Since $\|D^{2}v\|\leq\Delta v$, this proves
\begin{equation}\label{eqn5.2}
\|D^{2}v\|_{L^{\infty}((3/4)D)}
\leq C\big(1+\|\Delta v\|_{L^{p}(D)}^{p/n}\big),
\end{equation}
Here $\theta D$ denotes the concentric ball of relative radius $\theta$. The constant depends only on $n,k$, the radius of $D$, $m,M_{0}$ and $\omega$.

\begin{proof}[Proof of Theorem \ref{thm1.3}]
Fix the concentric balls in the theorem. Let $(u_i,f_i)$ be any sequence with common bounds for the stated data, including a common square-root Dini modulus $\omega_f$ on $B_{R}(x_{0})$. Take the subsequence and pairs $D_{\ell}\Subset B_{\ell}$ given by Lemma \ref{thm4.5}.

On each $B_{\ell}$, choose smooth Dirichlet approximations $U_{i,j}$ with $G(D^{2}U_{i,j})=g_{i,j}$. Obtain $g_{i,j}$ by extending $f_i^{1/k}$ by nearest-point projection onto $\overline B_\ell$ and mollifying. The projection is $1$-Lipschitz, so $g_{i,j}$ share the bounds $\lambda^{1/k},\Lambda^{1/k}$ and the modulus $(\lambda^{1/k-1}/k)\omega_f$, where $\lambda,\Lambda$ are the common bounds for $f_i$.

Lemma \ref{thm4.5} bounds $\|\Delta U_{i,j}\|_{L^{p}(D_\ell)}$ for $p=n+nk/2$, uniformly in $i$ and sufficiently large $j$ for each $i$. Estimate \eqref{eqn5.2} bounds the Hessians on $(3/4)D_\ell$, and interior interpolation with the supremum bound controls the gradients on $(1/2)D_\ell$. Passing to the limit in $j$ by uniform convergence and weak-$*$ compactness gives the same $W^{2,\infty}$ bounds for $u_i$. The half balls cover $\overline B_r(x_0)$.

A constant sequence gives local regularity. If the asserted estimate failed, a sequence with norms tending to infinity would contradict the subsequence bound.
\end{proof}

\subsection{H\"older densities}

\begin{proof}[Proof of Corollary \ref{thm1.4}]
By Section 3 and continuity, each $x_{0}\in\Omega$ is the center of two balls $D\Subset B\cap\{H_Bu>u\}$ and $B\Subset\Omega$, with $D$ of radius at most one. Use projection and mollification as above to obtain densities $g_j$ for the smooth Dirichlet approximants $U_j$. They share a positive lower bound $m$ and a modulus $Ct^\alpha$, which is square-root Dini. Theorem \ref{thm4.1} and \eqref{eqn5.2} thus bound their Hessians on $(3/4)D$ by $K_0$.

By translation and quadratic scaling take $(3/4)D=B_{1}$. For a rescaled approximant $v$ put $H=D^{2}v$ and $g=G(H)$.
The Hessians lie in the fixed compact convex set
\[
K=\{P\in\overline{\Gamma}_{k}:\ \|P\|\leq K_{0},\ G(P)\geq m\}\Subset\Gamma_{k}.
\]
By compactness and ellipticity, $\lambda_{*}I\leq DG\leq\Lambda_{*}I$ on $K$ for some positive $\lambda_{*},\Lambda_{*}$. Concavity and homogeneity give $G=\widehat G$ on $K$, where
\[
\widehat G(P)=\inf_{A\in K}DG(A)[P]
\]
is concave and uniformly elliptic on the space of symmetric matrices.
Subtract an affine function so that $\|v\|_{L^{\infty}(B_{1})}\leq CK_{0}$. The Evans-Krylov theorem \cite{Evans,Krylov} and Caffarelli's interior estimate \cite{CaffFully}, in the form of Caffarelli-Cabr\'e \cite[Theorems 6.6 and 8.1]{CC}, bound $[H]_{C^{0,\gamma}(B_{7/8})}$ uniformly for some $0<\gamma<\alpha$.

To recover $\alpha$, freeze $DG$ at $H(x)$, as in the author \cite[Section 3.4]{LyuSigma2Global} and \cite[Section 3.5]{LyuKGlobal}.
Put $D=B_{3/4}$, $d(x)=\operatorname{dist}(x,\partial D)$, and let
\[
M=\sup_{x\ne y\in D}\min\{d(x),d(y)\}^{\alpha}
\frac{\|H(x)-H(y)\|}{|x-y|^{\alpha}} .
\]
For $x\in D$ and $s=\eta d(x)$ with $0<\eta<1/4$, write $H_{0}=H(x)$. On $B_{s}(x)$,
\[
DG(H_{0})[H-H_{0}]=g-g(x)-E,\quad
E=G(H)-G(H_{0})-DG(H_{0})[H-H_{0}] .
\]
The bound for $D^{2}G$ on $K$ and the $C^{0,\gamma}$ estimate give $[E]_{C^{0,\alpha}(B_{s}(x))}\leq C\eta^{\gamma}[H]_{C^{0,\alpha}(B_{s}(x))}$.
After subtracting the second Taylor polynomial of $v$ at $x$, the Schauder estimate (see Gilbarg-Trudinger \cite{GT}) gives
\[
[H]_{C^{0,\alpha}(B_{s/2}(x))}
\leq C\big(s^{-\alpha}K_{0}+[g]_{C^{0,\alpha}}+
\eta^{\gamma}[H]_{C^{0,\alpha}(B_{s}(x))}\big).
\]
Order each pair so that $d(x)\leq d(y)$. For $|x-y|\geq\eta d(x)/2$, its contribution to $M$ is at most $2K_{0}(2/\eta)^{\alpha}$. For closer pairs the last estimate applies and $d\geq d(x)/2$ on $B_{s}(x)$. Thus
$M\leq C_{\eta}(K_{0}+[g]_{C^{0,\alpha}})+C\eta^{\gamma}M$.
Since $M$ is finite for each smooth approximant, choose $\eta$ small and absorb the last term to bound $[H]_{C^{0,\alpha}(B_{1/2})}$ uniformly.

The supremum and Hessian bounds control the gradients on smaller balls. Compactness in $C^{2}$ and uniform convergence to $u$ pass the $\alpha$-H\"older bound to $D^{2}u$.
\end{proof}

\section{Global $W^{2,1}$ regularity}

\subsection{Boundary extension}

We bound the Hessian mass by extending $u$ as a bounded $k$-convex function across $\partial\Omega$.

\begin{proof}[Proof of Theorem \ref{thm1.5}]
The geometric data in the theorem control all domain-dependent constants. Let $s$ be the signed distance, positive outside $\Omega$, and $N$ the outward unit normal. Following Caffarelli-Nirenberg-Spruck \cite[Section 2]{CNS} and the author \cite[Section 2.1]{LyuSigma2Global,LyuSigma2Higher,LyuKGlobal}, we build a defining function with Hessian uniformly in $\Gamma_{k}$ almost everywhere on a collar. The $C^{1,1}$ boundary charts give
\[
|(z-y)\cdot N(y)|\leq C|z-y|^{2},\quad \forall y,z\in\partial\Omega.
\]
For $x=y+tN(y)$ with $|t|$ small,
$|x-z|^{2}\geq t^{2}+(1-2C|t|)|z-y|^{2}$, so the nearest boundary point is $y$. These normal coordinates cover a collar. If $x'=z+t'N(z)$, then
\[
(x-x')\cdot(y-z)\geq\big(1-C(|t|+|t'|)\big)|y-z|^{2}.
\]
The coordinates are bi-Lipschitz locally, and $\nabla s(y+tN(y))=N(y)$ gives $s\in C^{1,1}$ on the collar. At almost every $y$, the tangential eigenvalues of $DN$ are $\kappa_{i}(y)$. Differentiating the distance formula yields the eigenvalues of $D^{2}s$ at almost every $y+tN(y)$,
\[
\frac{\kappa_{1}(y)}{1+t\kappa_{1}(y)},\ \ldots,\ 
\frac{\kappa_{n-1}(y)}{1+t\kappa_{n-1}(y)},\ 0.
\]
Indeed, the tangential differential of the coordinate map has eigenvalues $1+t\kappa_{i}(y)$, and the map preserves null sets.

The curvature vectors range almost everywhere in a compact subset of $\Gamma_{k-1}$. Since
\[
\sigma_{j}(\kappa,2J)=\sigma_{j}(\kappa)+2J\sigma_{j-1}(\kappa),\quad \forall 1\leq j\leq k,
\]
we may choose $J>0$ so that $(\kappa,2J)$ stays in a compact subset of $\Gamma_{k}$. The Hessian of $\rho:=s+Js^{2}$ has tangential eigenvalues $(1+2Jt)\kappa_{i}/(1+t\kappa_{i})$ and normal eigenvalue $2J$. Hence there are $\delta,\theta>0$, depending only on $k$ and $\Omega$, such that $J\delta\leq1/2$ and, almost everywhere on $\{|s|<\delta\}$,
\begin{equation}\label{eqn6.1}
D^{2}\rho-\theta I\in\Gamma_{k},\quad G(D^{2}\rho-\theta I)\geq\theta .
\end{equation}

Extend $\varphi$ in boundary charts and patch the extensions to obtain $\Phi\in C^{1,1}$ on the collar with $\|\Phi\|_{C^{1,1}}\leq C\|\varphi\|_{C^{1,1}(\partial\Omega)}$. Put $M:=\|\Phi\|_{C^{1,1}}+\Lambda^{1/k}$. Subharmonicity gives $u\leq\max_{\partial\Omega}\varphi$. For $\Omega\subset B_{R}(x_{0})$, comparison with
\[
\min_{\partial\Omega}\varphi+
\frac{\Lambda^{1/k}}{2\binom{n}{k}^{1/k}}\big(|x-x_{0}|^{2}-R^{2}\big),
\]
whose $k$-Hessian is $\Lambda$, gives the lower bound and hence $\|u\|_{L^{\infty}(\Omega)}\leq CM$ (see Trudinger-Wang \cite[Theorem 3.1]{TW1}).

Set $v:=\Phi+A\rho$, where $A=CM$. For large $C$ we have $D^{2}\Phi+A\theta I\geq0$ almost everywhere, so \eqref{eqn6.1} and monotonicity give $D^{2}v\in\Gamma_{k}$ and $G(D^{2}v)\geq A\theta\geq\Lambda^{1/k}$. By local mollification, concavity and weak continuity (see Trudinger-Wang \cite[Theorem 1.1]{TW1}), $v$ is $k$-convex and $\mu_{k}[v]\geq\Lambda\,dx$ on the collar. On $\partial\Omega$ we have $v=u$. On $\{s=-\delta/2\}$ we have $\rho\leq-\delta/4$, so $v\leq u$ for large $C$. Comparison gives $v\leq u$ on $\{-\delta/2<s<0\}$.

Define $U=u$ in $\Omega$ and $U=v$ on the rest of $\Omega^{+}:=\Omega\cup\{|s|<\delta/2\}$. This gives a continuous extension with $\|U\|_{L^{\infty}(\Omega^{+})}\leq CM$. A test function touching $U$ from above at $\partial\Omega$ also touches $v$, since $v\leq U$ near $\partial\Omega$ and $v=U$ there. Thus $U$ is $k$-convex on $\Omega^{+}$.

Lemma \ref{thm2.1} gives $U\in W_{\mathrm{loc}}^{1,1}(\Omega^{+})$ and $|D^{2}U|\leq\Delta U$ as measures. Take a cutoff $0\leq\eta\leq1$ in $C_{c}^{\infty}(\Omega^{+})$ equal to one near $\overline{\Omega}$, with derivatives controlled by the collar width. By Theorem \ref{thm1.1},
\[
\|D^{2}u\|_{L^{1}(\Omega)}
\leq\Delta U(\Omega)
\leq\int_{\Omega^{+}}\eta\,d\Delta U
=\int_{\Omega^{+}}U\Delta\eta\,dx\leq CM.
\]
Cover $\overline{\Omega}$ by balls of radius $\delta/16$ with centers in a maximal $\delta/16$-separated subset of $\overline{\Omega}$. Their doubles lie in $\Omega^{+}$ and the number of balls is controlled by $n$, $\delta$ and the diameter. The cutoff and weighted-mean Poincar\'e estimates of Section 2 apply on each ball to $U$, using $\nabla U\in BV_{\mathrm{loc}}(\Omega^{+})$ and $|D^{2}U|\leq\Delta U$. Thus $\|\nabla u\|_{L^{1}(\Omega)}\leq CM$ and
\[
\|u\|_{W^{2,1}(\Omega)}
\leq C\big(\|\varphi\|_{C^{1,1}(\partial\Omega)}+\Lambda^{1/k}\big).
\]
The positive lower bound for $f$ enters only through the interior absolute continuity theorem.
\end{proof}

\subsection{The boundary counterexample}

The boundary trace of a $2$-convex $W^{2,1}$ function has a spherical Laplacian with nonnegative singular part.

\begin{lemma}\label{thm6.1}
Let $u\in W^{2,1}(B_{1})\cap C^{0}(\overline{B}_{1})$ be $2$-convex. On $S=\partial B_{1}$ let $\varphi$ and $g\in L^{1}(S)$ be the traces of $u$ and its outward normal derivative, respectively. Then
\begin{equation}\label{eqn6.2}
\Delta_{S}\varphi+(n-1)g\geq0\quad\text{in distributions on }S.
\end{equation}
In particular, $\Delta_{S}\varphi$ is a signed measure and $(\Delta_{S}\varphi)^{-}\leq(n-1)g_{+}\,dS$.
\end{lemma}
\begin{proof}
For almost every $1/2<r<1$, the polar coordinate identity gives
\[
\Delta_{S}[u(r\,\cdot)]+(n-1)r(\partial_{r}u)(r\,\cdot)
=r^{2}(\Delta u-\partial_{rr}u)(r\,\cdot)\geq0
\]
in distributions on $S$, by smooth approximation on compact annuli and \eqref{eqn2.1}. For almost every $1/2<r<s<1$,
\[
\|(\partial_{r}u)(s\,\cdot)-(\partial_{r}u)(r\,\cdot)\|_{L^{1}(S)}
\leq\int_{r}^{s}\int_{S}\|D^{2}u(t\xi)\|\,dS(\xi)\,dt .
\]
As $r,s\uparrow1$, the right-hand side tends to zero. Thus the radial derivatives converge in $L^{1}(S)$ to $g$. Also $u(r\,\cdot)\to\varphi$ uniformly, giving \eqref{eqn6.2}. Its left-hand side is a nonnegative Radon measure, which yields the bound for the negative part.
\end{proof}

\begin{example}\label{thm6.2}
For each $2\leq k\leq n$ there is a boundary value $\varphi\in\bigcap_{0<\beta<1}C^{1,\beta}(\partial B_{1})$ for which the continuous $k$-convex solution of $\mu_{k}[u]=dx$ in $B_{1}$ lies outside $W^{2,1}(B_{1})$.
\end{example}
\begin{proof}
Starting with $[-1/2,1/2]$, at step $j\geq1$ retain two end intervals of length $\ell_{j}:=2^{-j}/(j+1)$ in each interval of step $j-1$. The probability measure $m$ giving mass $2^{-j}$ to each interval of step $j$ is atomless and supported on a set of length zero, since $2^{j}\ell_{j}\to0$. Let $F(t):=m((-\infty,t])$ and $\ell_{0}=1$. An interval of length $\ell_{j+1}<h\leq\ell_{j}$ meets at most two intervals of step $j$, and $2^{-j}\leq Ch\log(e/h)$. Hence
\[
|F(t)-F(s)|\leq C|t-s|\log(e/|t-s|),\quad \forall 0<|t-s|\leq1 .
\]
Its primitive $G_{0}(t):=\int_{-1/2}^{t}F(s)\,ds$ is locally $C^{1,\beta}$ for every $0<\beta<1$ and has distributional second derivative $m$.

Set $\varphi(\xi):=-G_{0}(\xi_{1})$. For a smooth function $b$ of $t=\xi_{1}$, the spherical Laplacian is $(1-t^{2})b''(t)-(n-1)tb'(t)$. In the band $|t|<3/4$ the coordinates $\xi=(t,\sqrt{1-t^{2}}\,\theta)$, $\theta\in S^{n-2}$, give $dS(\xi)=(1-t^{2})^{(n-3)/2}\,dt\,dS(\theta)$. For $n=2$, $S^{0}=\{-1,1\}$ carries counting measure. Mollifying $G_{0}$ in one variable now gives the singular part of $\Delta_{S}\varphi$ in the band,
\[
-(1-t^{2})^{(n-1)/2}\,dm(t)\,dS(\theta).
\]
This nonzero negative measure violates the trace condition in Lemma \ref{thm6.1}.

For smooth mollifications $G_{j}$ of $G_{0}$, let $u_{j}$ solve the Dirichlet problem with density one and boundary values $\varphi_{j}(\xi)=-G_{j}(\xi_{1})$ (see Caffarelli-Nirenberg-Spruck \cite{CNS}). By comparison (see Trudinger-Wang \cite[Theorem 3.1]{TW1}),
\[
\|u_{i}-u_{j}\|_{L^{\infty}(B_{1})}\leq\|\varphi_{i}-\varphi_{j}\|_{L^{\infty}(S)} .
\]
The uniform limit $u$ is $k$-convex with trace $\varphi$ and $\mu_{k}[u]=dx$ by weak continuity (see Trudinger-Wang \cite[Theorem 1.1]{TW1}). Comparison gives uniqueness. The trace obstruction applies because $k\geq2$.

For each $0<\beta<1$, the boundary values $\varphi_{j}$ are bounded in $C^{1,\beta}(S)$ and $\|u_{j}\|_{W^{2,1}(B_{1})}\to\infty$. Indeed, if a subsequence had bounded norms, the first and second derivatives of $u$ would have finite total variation in $B_{1}$. By Theorem \ref{thm1.1} they are absolutely continuous there. This would give $u\in W^{2,1}(B_{1})$.
\end{proof}

\appendix

\section{A Dirichlet problem with trace-dependent density}

We solve the comparison problem used in Section 4.3. The allowable change in density is independent of the weight and cutoff level.

Let $\eta\in C^{\infty}(\mathbb{R})$ be nondecreasing, equal to zero on $(-\infty,1]$ and one on $[2,\infty)$, and, for $m>0$, put
\[
C_{\eta}:=\max\left\{1,\sup_{z>0}(z^{2}\eta''(z))_{+}\right\},
\quad d_{*}:=\frac{m}{2kC_{\eta}} .
\]

\begin{proposition}\label{thmA.1}
Let $2\leq k\leq n$, $B=B_{R}(0)$, and let $U,g,b,\Theta\in C^{\infty}(\overline{B})$ satisfy
\[
\Theta>0,\quad D^{2}U\in\Gamma_{k},\quad G(D^{2}U)=g,\quad 0<m\leq b\leq g\leq M .
\]
If $\|g-b\|_{L^{\infty}}\leq d_{*}$, then for every $P>0$ the problem
\[
G(D^{2}v)=g\exp\left[-\log(g/b)\,
\eta\left(\frac{\Theta\Delta v}{P}\right)\right]
\quad\text{in }B,\qquad v=U\quad\text{on }\partial B
\]
has a unique admissible solution $v\in C^{2}(\overline B)\cap C^{\infty}(B)$ with $U\leq v$ and $m\leq G(D^{2}v)\leq M$.
\end{proposition}
\begin{proof}
We adapt the result of Chen-Zhou-Zhu \cite[Theorem 3.1]{CZZ}. See also the author \cite[Propositions 2.3 and 3.1]{LyuSigma2Interior}, \cite[Appendix A, (A.9)]{LyuSigma2Higher} and \cite[Section 2.3]{LyuKGlobal}. Put $\ell:=\log(g/b)$. For $0\leq t\leq1$ define
\[
h(x,T):=g(x)\exp[-t\ell(x)\eta(\Theta(x)T/P)],\quad
F(x,H):=\log\sigma_{k}(H)-k\log h(x,\operatorname{tr}H).
\]
Since $h>0$ and $h_{T}\leq0$, both $F$ and $G(H)-h(x,\operatorname{tr}H)$ are elliptic on $\Gamma_{k}$. Put $T=\operatorname{tr}H$, $z=\Theta T/P$ and $Q=(\sigma_{k}/\sigma_{1})^{1/(k-1)}$. Quotient concavity gives concavity of $\log Q$, and
\[
F=(k-1)\log Q+\log T-k\log g+kt\ell\eta(z).
\]
Since $0\leq\ell\leq\|g-b\|_{L^{\infty}}/m$, differentiation at fixed $x$ yields
\[
D_{H}^{2}F(x,H)[Z,Z]
\leq[-1+kt\ell z^{2}\eta''(z)]\left(\frac{\operatorname{tr}Z}{T}\right)^{2}
\leq-\frac12\left(\frac{\operatorname{tr}Z}{T}\right)^{2}.
\]

Let $v\in C^{2}(\overline B)\cap C^{\infty}(B)$ be admissible with $F(x,D^{2}v)=0$ and boundary value $U$. In the estimates below the constants may depend on $P$ and the fixed smooth data, but are uniform in $t$. Let $V$ solve $G(D^{2}V)=m/2$ with boundary value $U$. Comparison gives $U\leq v\leq V$. On $\partial B$, the tangential gradient is prescribed and $V_{\nu}\leq v_{\nu}\leq U_{\nu}$ for the outward unit normal $\nu$.

Write $L=A^{ij}\partial_{ij}$, where $A=DG(D^{2}v)-h_{T}(x,\Delta v)I$. Homogeneity and concavity give $\operatorname{tr}A\geq DG(D^{2}v)[I]\geq G(I)$. At fixed $T$, $D_{x}h$ is bounded uniformly in $T$ and $t$, since $\Theta T/P\leq2$ on the support of $\eta'$. For a rotation field $X=x_{i}\partial_{j}-x_{j}\partial_{i}$, rotational invariance of $G$ and the trace gives $L(Xv)=X_{x}h(x,\Delta v)$, with the trace held fixed on the right. If $w=v-U$, then $Xw$ vanishes on $\partial B$ and $|L(Xw)|\leq C\operatorname{tr}A$. Comparison with $\pm C(R^{2}-|x|^{2})$ gives $|Xw|\leq C(R^{2}-|x|^{2})$. On the sphere $\nabla w=w_{\nu}\nu$ and $D_{\nu}X=X/R$, whence
\[
\partial_{\nu}(Xw)=D^{2}w(\nu,X)+(X/R)\cdot\nabla w=D^{2}w(\nu,X).
\]
This bounds the mixed entries of $D^{2}v$ on $\partial B$, since the rotations span every tangent space there.

The boundary data and the bound for $v_{\nu}$ control the tangential block $H'=D^{2}v|_{\nu^{\perp}}$. Since $v\leq V$ with equality on $\partial B$,
\[
H'-D^{2}V|_{\nu^{\perp}}=R^{-1}(v-V)_{\nu}I\geq0 .
\]
The tangential block of $D^{2}V$ stays in a compact subset of $\Gamma_{k-1}$, so ellipticity gives $\sigma_{k-1}(H')\geq c>0$. This is the coefficient of $v_{\nu\nu}$ in $\sigma_{k}(D^{2}v)=h^{k}\leq M^{k}$. All other terms involve only bounded entries, so $v_{\nu\nu}$ is bounded on $\partial B$.

To pass the boundary bound to the interior, put $q_{t}:=g^{1-t}b^{t}$ and choose $K>0$ so that
\[
2KG(I)>\sup_{0\leq t\leq1}\|\Delta q_{t}\|_{L^{\infty}(B)} .
\]
At an interior maximum of $T+K|x|^{2}$ with $T=\Delta v>2P/\min_{\overline B}\Theta$, the equation is $G(D^{2}v)=q_{t}$ in a neighborhood of the point. Taking the Laplacian and using concavity gives
\[
\Delta q_{t}\leq DG(D^{2}v)[D^{2}T]\leq-2K G(I),
\]
a contradiction. Thus $\Delta v$ and hence $\|D^{2}v\|$ are bounded throughout $\overline B$. Integrating from the boundary gives a uniform $C^{2}(\overline B)$ bound. Since $G(D^{2}v)\geq m$, the Hessians range in a compact subset of $\Gamma_{k}$ where $F$ has uniform ellipticity and derivative bounds through order two.

To apply the global concave estimate, start with a $C^{2,\gamma}(\overline B)$ solution, $0<\gamma<1$. Interior regularity (see Gilbarg-Trudinger \cite[Lemma 17.16]{GT}) gives smoothness in $B$. For each rotation field $X$, the equation $L(Xv)=X_{x}h(x,\Delta v)$ has $C^{0,\gamma}$ coefficients and right-hand side on $\overline B$ and boundary value $XU$. The linear Dirichlet theorem and uniqueness then give $Xv\in C^{2,\gamma}(\overline B)$. Thus in a smooth tangential-radial frame near $\partial B$, every Hessian entry with a tangential index is $C^{1,\gamma}$. Since $F_{H}[\nu\otimes\nu]>0$, the equation gives the same regularity for $v_{\nu\nu}$, and $v\in C^{3,\gamma}(\overline B)$.

The global estimate (see Gilbarg-Trudinger \cite[Theorem $17.26'$]{GT}) gives a uniform $C^{2,\alpha}(\overline B)$ bound for some $\alpha>0$. Use continuity in $C^{2,\beta}(\overline B)$ with $0<\beta<\alpha$, starting at $U$ for $t=0$. The linearization has no zeroth-order term, so the linear Dirichlet isomorphism gives openness. Compactness and the fixed compact Hessian range in $\Gamma_{k}$ give closedness and hence existence at $t=1$. Linearizing along the segment between the Hessians of two solutions gives uniqueness by the maximum principle.
\end{proof}

\section*{Acknowledgements}

The author is deeply grateful to Professors Dongsheng Li and Yu Yuan for their guidance and support.
AI tools were used as an aid in preliminary mathematical exploration and in refining the exposition.
All mathematical arguments were independently verified by the author, who takes full responsibility for the manuscript.

\end{document}